\documentclass[leqno,draft]{article}

\newtheorem{theorem}{Theorem}
\newtheorem{lemma}[theorem]{Lemma}
\newtheorem{proposition}[theorem]{Proposition}
\newtheorem{definition}[theorem]{Definition}
\newtheorem{corollary}[theorem]{Corollary}

\newcommand{\begintheorem}{\addtocounter{equation}{1}\begin{theorem}}
\newcommand{\beginlemma}{\addtocounter{equation}{1}\begin{lemma}}
\newcommand{\beginproposition}{\addtocounter{equation}{1}\begin{proposition}}
\newcommand{\begindefinition}{\addtocounter{equation}{1}\begin{definition}}
\newcommand{\begincorollary}{\addtocounter{equation}{1}\begin{corollary}}

\begin{document}

\title{Some aspects of algebraic topology: \\ An analyst's
perspective}

\author{Stephen Semmes	\\
	Rice University}

\date{}

\maketitle

\begin{abstract}
These informal notes discuss a few basic notions and examples, with
emphasis on constructions that may be relevant for analysis on metric
spaces.
\end{abstract}

\tableofcontents

\bigskip

	From the point of view of classical analysis, the story of
algebraic topology is fantastic.  In particular, the theory of
homology and cohomology is a far-reaching extension of calculus.  Here
the focus is on slightly different but closely related objects.

\section{Open coverings}
\label{open coverings}
\setcounter{equation}{0}

	Let $X$ be a topological space, and let $\mathcal{U} =
\{U_\alpha\}_{\alpha \in A}$ be a collection of open subsets of $X$.
We say that $\mathcal{U}$ is an \emph{open covering} of a set $E
\subseteq X$ if
\begin{equation}
	E \subseteq \bigcup_{\alpha \in A} U_\alpha.
\end{equation}
In particular, $\mathcal{U}$ is an open covering of $X$ if
\begin{equation}
	\bigcup_{\alpha \in A} U_\alpha = X.
\end{equation}

	One can think of an open covering as a kind of geometric
structure, which specifies a level of localization at each point in
the set.  In a metric space, one might consider coverings by balls of
some radii, or conditions on the diameters of the open sets in a
covering.

	As usual, a set $E \subseteq X$ is said to be \emph{compact}
if for every open covering $\{U_\alpha\}_{\alpha \in A}$ of $E$ in $X$
there are finitely many indices $\alpha_1, \ldots, \alpha_n \in A$
such that
\begin{equation}
	E \subseteq U_{\alpha_1} \cup \cdots \cup U_{\alpha_n}.
\end{equation}
Similarly, $E$ is \emph{countably compact} if for each open covering
$\{U_\alpha\}_{\alpha \in A}$ of $E$ there is an $A_1 \subseteq A$
with only finitely or countably many elements such that
\begin{equation}
	E \subseteq \bigcup_{\alpha \in A_1} U_\alpha.
\end{equation}
Thus compact sets are automatically countably compact.

	It is well known that a metric space $M$ is countably compact
if and only if it is \emph{separable}, which is to say that there is a
dense set in $M$ with only finitely or countably many elements.  These
properties are also equivalent to the existence of a base for the
topology of $M$ with only finitely or countably many elements.  The
existence of a base for the topology of any topological space with
only finitely or countably many elements implies that the space is
countably compact and contains a dense set with only finitely or
countably many elements.

	For example, the set of rationals is a countable dense set in
the real line ${\bf R}$ with the standard topology, and the collection
of open intervals with rational endpoints is a countable base for the
topology.  The compact subsets of the real line are exactly the closed
and bounded sets.

	Suppose that $\mathcal{U} = \{U_\alpha\}_{\alpha \in A}$,
$\mathcal{V} = \{V_\beta\}_{\beta \in B}$ are families of open subsets
of a topological space $X$.  We say that $\mathcal{V}$ is a
\emph{refinement} of $\mathcal{U}$ if for every $\beta \in B$ there is
an $\alpha \in A$ such that $V_\beta \subseteq U_\alpha$.  We may
express this relationship with the notation $\mathcal{U} \prec
\mathcal{V}$, where the ordering reflects the fact that a refinement
$\mathcal{V}$ of $\mathcal{U}$ corresponds to greater precision in
$X$.

	Let $\mathcal{U}$, $\mathcal{V}$, and $\mathcal{W}$ be
families of open subsets of $X$.  If $\mathcal{U} \prec \mathcal{V}$
and $\mathcal{V} \prec \mathcal{W}$, then it is easy to see that
$\mathcal{U} \prec \mathcal{W}$.

	Let $\mathcal{U} = \{U_\alpha\}_{\alpha \in A}$ and
$\mathcal{V} = \{V_\beta\}_{\beta \in B}$ be families of open subsets
of $X$, and consider the collection $\mathcal{W}$ of open subsets of
$X$ of the form $U_\alpha \cap V_\beta$ for $\alpha \in A$ and $\beta
\in B$.  Clearly
\begin{equation}
	\mathcal{U}, \mathcal{V} \prec \mathcal{W},
\end{equation}
and any set $E \subseteq X$ which is covered by $\mathcal{U}$ and by
$\mathcal{V}$ is also covered by $\mathcal{W}$.

\section{Connectedness}
\label{connectedness}
\setcounter{equation}{0}

	Let $X$ be a topological space, and suppose that $\mathcal{U}
= \{U_\alpha\}_{\alpha \in A}$ is a collection of open subsets of $X$.

	Let us say that $\alpha, \alpha' \in A$ are \emph{adjacent} if
\begin{equation}
	U_\alpha \cap U_{\alpha'} \ne \emptyset.
\end{equation}
An equivalence relation $\sim$ can be defined on $A$ by putting
$\alpha \sim \beta$ when $\alpha, \beta \in A$ and there is a finite
sequence $\alpha_1, \ldots, \alpha_n$ of elements of $A$ such that
$\alpha_1 = \alpha$, $\alpha_n = \beta$, and $\alpha_i$ is adjacent to
$\alpha_{i + 1}$ for $1 \le i < n$.

	If $E \subseteq X$ is connected and $\mathcal{U}$ is an open
covering of $E$, then $\alpha \sim \beta$ for every $\alpha, \beta \in
A$ such that
\begin{equation}
	E \cap U_\alpha \ne \emptyset
		\quad\hbox{and}\quad E \cap U_\beta \ne \emptyset.
\end{equation}
For let $\alpha \in A$ with $E \cap U_\alpha \ne \emptyset$ be given,
let $V$ be the union of the $U_\beta$'s with $\beta \in A$ and $\alpha
\sim \beta$, and let $W$ be the union of the $U_\beta$'s with $\beta
\in A$ and $\alpha \not\sim \beta$.  By construction, $V$ and $W$ are
disjoint open subsets of $X$.  Also, $E \cap V \ne \emptyset$, since
$U_\alpha \subseteq V$.  The connectedness of $E$ implies that $E \cap
W = \emptyset$, and hence $E \cap U_\beta = \emptyset$ when $\alpha
\not\sim \beta$.

	Conversely, if $E$ is a disconnected open set in $X$, then $E$
can be expressed as the union of two nonempty disjoint open sets.  The
open covering of $E$ consisting of these two open sets does not have
the property just described.

	Similarly, if $E$ is a disconnected closed set in $X$, then
$E$ can be expressed as the union of two nonempty disjoint closed
sets.  If $X$ is normal, then these two disjoint closed sets are
contained in a pair of disjoint open sets which cover $E$ and do not
have the property just described.

	Suppose now that $X$ is \emph{Hausdorff}, which means that for
every $p, q \in X$ with $p \ne q$ there are disjoint open subsets $V$,
$W$ of $X$ such that $p \in V$ and $q \in W$.  Let $K$ be a compact
set in $X$, and let $p$ be an element of $X$ not in $K$.  For each $q
\in K$, let $V(q)$, $W(q)$ be disjoint open subsets of $X$ such that
$p \in V(q)$ and $q \in W(q)$.  By compactness, there are finitely
many elements $q_1, \ldots, q_n$ of $K$ such that $K \subseteq
\bigcup_{i = 1}^n W(q_i)$.  If $V = \bigcap_{i = 1}^n V(q_i)$ and $W =
\bigcup_{i = 1}^n W(q_i)$, then $V$, $W$ are disjoint open subsets of
$X$ such that $p \in V$ and $K \subseteq W$.  In particular, it
follows that compact subsets of Hausdorff topological spaces are
closed.  In the same way, one can show that disjoint compact subsets
of a Hausdorff topological space are contained in disjoint open sets.

	If $K$ is a disconnected compact set in $X$, then $K$ can be
expressed as the union of two nonempty disjoint compact sets.  These
disjoint compact sets are contained in disjoint open sets, which form
an open covering of $K$ that does not have the property described
previously.

	In spaces like metric spaces, in which separated sets are
contained in disjoint open sets, similar remarks apply to arbitrary
disconnected sets.

\section{Topological dimension}
\label{topological dimension}
\setcounter{equation}{0}

	Let $X$ be a topological space, let $\mathcal{U} =
\{U_\alpha\}_{\alpha \in A}$ be a collection of open subsets of $X$,
and let $n$ be a nonnegative integer.  For each $x \in X$, put
$A(x) = \{\alpha \in A : x \in U_\alpha\}$.

	We say that $\mathcal{U}$ has \emph{order $\le n$} if $A(x)$
has $\le n + 1$ elements for every $x \in X$.  Equivalently,
\begin{equation}
	U_{\alpha_1} \cap \cdots \cap U_{\alpha_{n + 2}} = \emptyset
\end{equation}
for any $n + 2$ distinct elements $\alpha_1, \ldots, \alpha_{n + 2}$
of $A$.

	Roughly speaking, topological dimension $\le n$ corresponds to
coverings of order $\le n$ and arbitrarily high precision.  However,
one ought to be careful about some technical details.  This can
include additional hypotheses on the sets whose topological dimension
is being considered, as well as restrictions to coverings by finitely
many open sets, for instance.

	Note that a collection of open sets has order $0$ if and only
if the open sets in the collection are pairwise disjoint.  A set with
topological dimension $0$ is totally disconnected, which is to say
that the only connected subsets that it contains are the empty set and
sets with one element.

	It is easy to make coverings of subsets of the real line by
families of small open intervals such that no point is contained in
more than two of the intervals.  As a basic exercise, one can check
that every open covering of the closed unit interval $[0, 1]$ has a
refinement of order $1$ that is also a covering of $[0, 1]$.

	In the plane, one can certainly make coverings by small
neighborhoods of small squares.  If one uses squares in a standard
grid, then the vertices of the squares will be contained in four
adjacent squares.  One can avoid the problem by shifting some of the
squares so that no point is contained in more than three open sets.

	Alternatively, one can use the geometry of the standard grid
but choose the coverings a bit differently.  One can cover the
vertices separately, then the interiors of the edges, and then the
interiors of the squares.

\section{Continuous mappings}
\label{continuous mappings}
\setcounter{equation}{0}

	Let $X$, $Y$ be topological spaces, and let $f$ be a mapping
from $X$ to $Y$.  We say that $f$ is \emph{continuous} at a point $p
\in X$ if for every open set $V \subseteq Y$ with $f(p) \in V$ there
is an open set $U \subseteq X$ such that $p \in U$ and $f(x) \in V$
for every $x \in U$.

	We say that $f$ is a continuous mapping from $X$ to $Y$ if $f$
is continuous at every $p \in X$.  This is equivalent to the
requirement that for every open set $W \subseteq Y$,
\begin{equation}
	f^{-1}(W) = \{x \in X : f(x) \in W\}
\end{equation}
is an open set in $X$.

	For any mapping $f : X \to Y$ and set $E \subseteq Y$, $X
\backslash f^{-1}(E) = f^{-1}(Y \backslash E)$.  Here $A \backslash B$
means the collection of elements of a set $A$ which are not in a set
$B$.  It follows that $f : X \to Y$ is continuous if and only if
$f^{-1}(E)$ is a closed set in $X$ for every closed set $E$ in $Y$.

	The space of continuous mappings from $X$ to $Y$ is denoted
$\mathcal{C}(X, Y)$.  This may be abbreviated to $\mathcal{C}(X)$ when
$Y$ is the real line with the standard topology.

	If $X$, $Y$, $Z$ are topological spaces and $f : X \to Y$, $g
: Y \to Z$ are continuous mappings, then the composition $g \circ f$
is defined by $(g \circ f)(x) = g(f(x))$ and is a continuous mapping
from $X$ to $Z$.  This follows from any of the aforementioned
characterizations of continuity.

	A continuous mapping $f : X \to Y$ is a \emph{homeomorphism}
from $X$ onto $Y$ if there is a continuous inverse mapping $g : Y \to
X$, which is to say that $g(f(x)) = x$ for every $x \in X$ and
$f(g(y)) = y$ for every $y \in Y$.

	If $E \subseteq X$ is connected and $f : X \to Y$ is continuous,
then $f(E)$ is connected in $Y$.

	Suppose that $f : X \to Y$ is continuous and $\mathcal{U} =
\{U_\alpha\}_{\alpha \in A}$ is a collection of open subsets of $Y$.
Let $f^*(\mathcal{U})$ be the collection of open subsets of $X$ of the
form $f^{-1}(U_\alpha)$, $\alpha \in A$.

	If $\mathcal{U}$ has order $\le n$ in $Y$, then
$f^*(\mathcal{U})$ has order $\le n$ in $X$.  If $\mathcal{V}$ is
another collection of open subsets of $Y$ which is a refinement of
$\mathcal{U}$, then $f^*(\mathcal{V})$ is a refinement of
$f^*(\mathcal{U})$ in $X$.

	For $E \subseteq X$, $f(E)$ is the set of $y \in Y$ of the
form $y = f(x)$ for some $x \in E$.  If $\mathcal{U}$ is an open
covering of $f(E)$ in $Y$, then $f^*(\mathcal{U})$ is an open covering
of $E$ in $X$.  One can use this to show that compactness or countable
compactness of $E$ in $X$ implies the same property of $f(E)$ in $Y$.

	If $X$, $Y$, $Z$ are topological spaces, $f : X \to Y$ and $g
: Y \to Z$ are continuous mappings, and $\mathcal{W}$ is a collection
of open subsets of $Z$, then it is easy to see that $(g \circ
f)^*(\mathcal{W}) = f^*(g^*(\mathcal{W}))$.

\section{Uniform continuity}
\label{uniform continuity}
\setcounter{equation}{0}

	Let $(M, d(x, y))$ and $(N, \rho(u, v))$ be metric spaces.  A
mapping $f$ from $M$ to $N$ is \emph{uniformly continuous} if for every
$\epsilon > 0$ there is a $\delta > 0$ such that
\begin{equation}
	\rho(f(x), f(y)) < \epsilon
\end{equation}
when $x, y \in M$ and $d(x, y) < \delta$.  Every continuous mapping $f
: M \to N$ is uniformly continuous when $M$ is compact, by a
well-known theorem.

	Uniform continuity of mappings between general topological
spaces is not defined without additional information.  However, there
are special cases in which some of the same features are present.

	Let $(N, \rho(u, v))$ be a metric space again, and let $X$,
$Y$ be topological spaces.  We shall be interested in continuous
mappings $f : X \times Y \to N$, where $X \times Y$ is equipped with
the product topology.

	Let $\epsilon > 0$ be given.  As a partial version of uniform
continuity at a point $y \in Y$, consider the condition that there be
an open set $V$ in $Y$ such that $y \in V$ and
\begin{equation}
\label{rho(f(x, y), f(x, z)) < epsilon}
	\rho(f(x, y), f(x, z)) < \epsilon
\end{equation}
for each $x \in X$ and $z \in V$.

	If $X$ is compact, then every continuous mapping $f : X \times
Y \to N$ has this property.  Indeed, for each $w \in X$ there are open
sets $U(w) \subseteq X$ and $V(w) \subseteq Y$ such that $w \in U(w)$,
$y \in V(w)$, and
\begin{equation}
	\rho(f(x, z), f(w, y)) < \epsilon / 2
\end{equation}
when $x \in U(w)$ and $z \in V(w)$.  This follows from the continuity
of $f$ at $(w, y)$ and the definition of the product topology.  Hence
\begin{eqnarray}
	\rho(f(x, y), f(x, z))
	  & \le & \rho(f(x, y), f(w, y)) + \rho(f(w, y), f(x, z ))	\\
	  & < & \epsilon / 2 + \epsilon / 2 = \epsilon		\nonumber
\end{eqnarray}
when $x \in U(w)$ and $z \in V(w)$, by applying the previous statement
to $(x, y)$ and to $(x, z)$.  Because $X$ is compact, there are
finitely many elements $w_1, \ldots, w_n$ of $X$ such that $X =
\bigcup_{i = 1}^n U(w_i)$, and therefore (\ref{rho(f(x, y), f(x, z)) <
epsilon}) holds with $V = \bigcap_{i = 1}^n V(w_i)$.

	Suppose that $Y$ is also a metric space.  In this event we can
reformulate the partial version of uniform continuity as saying that
for each $\epsilon > 0$ and $y \in Y$ there is $\delta(y) > 0$ such
that (\ref{rho(f(x, y), f(x, z)) < epsilon}) holds when $x \in X$, $z
\in Y$, and the distance from $y$ to $z$ is less than $\delta(y)$.

	If $Y$ is compact too, then for each $\epsilon > 0$ there is a
$\delta > 0$ such that (\ref{rho(f(x, y), f(x, z)) < epsilon}) holds
when $x \in X$, $y, z \in Y$, and the distance from $y$ to $z$ in $Y$
is less than $\delta$.  As in the previous paragraph, for each $y \in
Y$ there is a $\delta_2(y) > 0$ such that
\begin{equation}
	\rho(f(x, y), f(x, z)) < \epsilon / 2
\end{equation}
when $x \in X$, $z \in Y$, and the distance from $y$ to $z$ is less
than $\delta_2(y)$.  Let $B(y)$ be the open ball in $Y$ centered at
$y$ and with radius $\delta_2(y) / 2$.  By compactness, there are
finitely many elements $y_1, \ldots, y_r$ of $Y$ such that
\begin{equation}
\label{Y = bigcup_{i = 1}^r B(y_i)}
	Y = \bigcup_{i = 1}^r B(y_i).
\end{equation}
Put
\begin{equation}
	\delta = \min \{\delta_2(y_i) / 2 : 1 \le i \le r\},
\end{equation}
and suppose that $y, z \in Y$ and the distance from $y$ to $z$ is less
than $\delta$.  By (\ref{Y = bigcup_{i = 1}^r B(y_i)}), there is an
$i$ such that $1 \le i \le r$ and $y \in B(y_i)$.  Hence the distance
from $z$ to $y_i$ is less than or equal to the sum of the distances
from $z$ to $y$ and from $y$ to $y_i$, which is less than $\delta +
\delta_2(y_i) / 2 \le \delta_2(y_i)$. Thus
\begin{eqnarray}
	\rho(f(x, y), f(x, z))
	  & \le & \rho(f(x, y), f(x, y_i)) + \rho(f(x, y_i), f(x, z))	\\
	  & < & \epsilon / 2 + \epsilon / 2 = \epsilon	\nonumber
\end{eqnarray}
for every $x \in X$, as desired.

\section{The supremum metric}
\label{supremum metric}
\setcounter{equation}{0}

	Let $X$ be a topological space, and let $(N, \rho(u, v))$ be a
metric space.  A mapping $f : X \to N$ is said to be \emph{bounded} if
$f(X)$ is a bounded set in $N$.  The space of bounded continuous
mappings from $X$ to $N$ is denoted $\mathcal{C}_b(X, N)$.  The space
of bounded continuous real-valued functions on $X$ may be denoted as
$\mathcal{C}_b(X)$.  If $X$ is compact and $f : X \to N$ is
continuous, then $f(X)$ is a compact set in $N$ and hence bounded.
Thus $\mathcal{C}_b(X, N) = \mathcal{C}(X, N)$ when $X$ is compact.
For $f_1, f_2 \in \mathcal{C}_b(X, N)$, consider
\begin{equation}
	d_*(f_1, f_2) = \sup \{\rho(f_1(x), f_2(x)) : x \in X\}.
\end{equation}
The supremum makes sense here, because $f_1$ and $f_2$ are bounded
functions on $X$, which implies that $\rho(f_1, f_2)$ is bounded.  One
can check that this defines a metric on $\mathcal{C}_b(X, N)$, known
as the supremum metric.

\section{Norms on vector spaces}
\label{norms on vector spaces}
\setcounter{equation}{0}

	The \emph{absolute value} of a real number $r$ is denoted
$|r|$ and defined to be $r$ when $r \ge 0$ and $-r$ when $r \le 0$.
For every $r, t \in {\bf R}$,
\begin{equation}
	|r + t| \le |r| + |t|
\end{equation}
and
\begin{equation}
	|r \, t| = |r| \, |t|.
\end{equation}

	Let $V$ be a vector space over the real numbers.  A
\emph{norm} on $V$ is a real-valued function $\|v\|$ on $V$ such that
$\|v\| \ge 0$ for every $v \in V$, $\|v\| = 0$ if and only if $v = 0$,
\begin{equation}
	\|r \, v\| = |r| \, \|v\|
\end{equation}
for every $r \in {\bf R}$ and $v \in V$, and
\begin{equation}
	\|v + w\| \le \|v\| + \|w\|
\end{equation}
for every $v, w \in V$.

	For example, the absolute value function defines a norm on
${\bf R}$.  If $n$ is a positive real number, then the space ${\bf
R}^n$ of $n$-tuples $x = (x_1, \ldots, x_n)$ of real numbers is a
vector space with respect to coordinatewise addition and scalar
multiplication.  The standard Euclidean norm on ${\bf R}^n$ is defined
by
\begin{equation}
	\|x\|_2 = \Big(\sum_{i = 1}^n x_i^2 \Big)^{1/2}
\end{equation}
and satisfies the conditions described in the previous paragraph.
Moreover,
\begin{equation}
	\|x\|_p = \Big(\sum_{i = 1}^n |x_i|^p \Big)^{1/p}
\end{equation}
defines a norm on ${\bf R}^n$ for each real number $p$ with $1 \le p <
\infty$.  This extends to $p = \infty$ by putting
\begin{equation}
	\|x\|_\infty = \max(|x_1|, \ldots, |x_n|).
\end{equation}

	If $V$ is a real vector space and $\|v\|$ is a norm on $V$, then
\begin{equation}
	d(v, w) = \|v - w\|
\end{equation}
defines a metric on $V$.  Two norms $\|v\|$, $\|v\|'$ on a vector
space $V$ are said to be \emph{equivalent} if there is a positive real
number $C$ such that
\begin{equation}
	C^{-1} \, \|v\| \le \|v\|' \le C \, \|v\|
\end{equation}
for every $v \in V$.  In this event the corresponding metrics satisfy
the analogous relation, and they determine the same topology on $V$.
Conversely, two norms on $V$ are equivalent if the corresponding
metrics determine the same topology on $V$.

	The standard metric on the real line corresponds to the
absolute value function in this way, and the standard Euclidean metric
on ${\bf R}^n$ corresponds to the Euclidean norm $\|x\|_2$.  For each
positive integer $n$, the norms $\|x\|_p$, $1 \le p \le \infty$, are
equivalent on ${\bf R}^n$, with a constant that depends on $n$.  There
are similar norms on infinite-dimensional spaces that are not
equivalent.

	Let $X$ be a topological space, and consider the vector space
$\mathcal{C}_b(X)$ of bounded continuous real-valued functions on $X$.
The \emph{supremum norm} on $\mathcal{C}_b(X)$ is defined by
\begin{equation}
	\|f\|_* = \sup \{|f(x)| : x \in X\}.
\end{equation}
It is easy to check that $\|f\|_*$ is a norm on $\mathcal{C}_b(X)$,
for which the corresponding metric is the supremum metric, and that
\begin{equation}
	\|f_1 \, f_2\|_* \le \|f_1\|_* \, \|f_2\|_*
\end{equation}
for every $f_1, f_2 \in \mathcal{C}_b(X)$.  If $V$ is a real vector
space equipped with a norm $\|v\|$, then $\mathcal{C}(X, V)$ is a
vector space with respect to pointwise addition and scalar
multiplication, and $\mathcal{C}_b(X, V)$ is a linear subspace of
$\mathcal{C}(X, V)$.  Furthermore,
\begin{equation}
	\|f\|_{*, V} = \sup \{\|f(x)\| : x \in X\}
\end{equation}
defines a norm on $\mathcal{C}_b(X, V)$, and
\begin{equation}
	d_*(f_1, f_2) = \|f_1 - f_2\|_{*, V}
\end{equation}
is the corresponding supremum metric on $\mathcal{C}_b(X, V)$.

\section{Pathwise connected sets}
\label{pathwise connected sets}
\setcounter{equation}{0}

	Let $X$ be a topological space, and let $E$ be a set contained
in $X$.  We say that $E$ is \emph{pathwise connected} if for every $p,
q \in E$ there is a continuous mapping $\phi$ from a closed interval
$[a, b]$ in the real line into $X$ such that $\phi(a) = p$, $\phi(b) =
q$, and $\phi([a, b]) \subseteq E$.

	Intervals in the real line are connected sets, and continuous
mappings send connected sets to connected sets.  Hence the sets
$\phi([a, b])$ as in the previous paragraph are connected.  One can
use this to show that pathwise connected sets are automatically
connected.

	If $E \subseteq X$ is pathwise connected, $Y$ is another
topological space, and $f : X \to Y$ is continuous, then $f(E)$ is
pathwise connected in $Y$.  This follows from the fact that
compositions of continuous mappings are continuous.

	For any $E \subseteq X$, let $\sim_E$ be the relation on $E$
defined by $p \sim_E q$ when $p, q \in E$ and there is a continuous
path $\phi : [a, b] \to X$ such that $\phi(a) = p$, $\phi(b) = q$, and
$\phi([a, b]) \subseteq E$.  It is easy to see that $\sim_E$ is an
equivalence relation on $E$.  The equivalence classes in $E$
associated to $\sim_E$ are pathwise connected, and are known as the
pathwise connected components of $E$.

	Let $V$ be a real vector space.  A set $E \subseteq V$ is said
to be \emph{convex} if for every $v, w \in E$ and real number $t$
with $0 < t < 1$,
\begin{equation}
\label{t v + (1 - t) w}
	t \, v + (1 - t) \, w
\end{equation}
is also an element of $E$.

	Suppose that $V$ is equipped with a norm $\|v\|$, and hence a
metric $\|v - w\|$ and a topology.  Using the triangle inequality, one
can check that open and closed balls in $V$ are convex sets.

	If $E \subseteq V$ is convex, then $E$ is pathwise connected,
because (\ref{t v + (1 - t) w}) is a continuous function of $t \in
{\bf R}$ with values in $V$.

	Let $U$ be an open set in $V$, and let $\sim_U$ be the
equivalence relation on $U$ associated to continuous paths in $U$ as
before.  If $p \in U$, then there is an open ball $B$ in $V$ centered
at $p$ such that $B \subseteq U$, and $p \sim_U q$ for every $q \in
B$.  This implies that the pathwise connected components of $U$ are
open subsets of $V$.  By construction, the pathwise connected
components of any set are pairwise disjoint.  It follows that $U$ is
pathwise connected if it is connected.  If $V$ is separable, then any
collection of pairwise disjoint open subsets of $V$ has only finitely
or countably many elements.  Thus an open set $U \subseteq V$ has only
finitely or countably many pathwise connected components when $V$ is
separable.

	Of course, ${\bf R}^n$ is separable.  If $M$, $N$ are metric
spaces with $M$ compact and $N$ separable, then $\mathcal{C}(M, N)$ is
separable with respect to the supremum metric.  In particular, the
vector space $\mathcal{C}(M)$ is separable with respect to the
supremum norm when $M$ is a compact metric space.

\section{Homotopies}
\label{homotopies}
\setcounter{equation}{0}

	Let $X$, $Y$ be topological spaces, and let $f, g : X \to Y$
be continuous mappings.  A \emph{homotopy} between $f$ and $g$ is a
continuous mapping
\begin{equation}
	H : X \times [a, b] \to Y
\end{equation}
such that $H(x, a) = f(x)$ and $H(x, b) = g(x)$ for every $x \in X$.

	Thus a homotopy is basically a continuous path in
$\mathcal{C}(X, Y)$.  If $X$ is compact and $Y$ is a metric space,
then a homotopy is exactly a continuous path in $\mathcal{C}(X, Y)$
with respect to the supremum metric.

	We say that $f, g : X \to Y$ are \emph{homotopic} if there is
a homotopy between them.  One can check that homotopy defines an
equivalence relation on $\mathcal{C}(X, Y)$.  The corresponding
equivalence classes are known as \emph{homotopy classes}.

	If $X$ has only one element, then a mapping from $X$ to $Y$ is
essentially the same as a point in $Y$.  Homotopy classes of these
mappings correspond to the pathwise connected components of $Y$.

	Suppose that $V$ is a real vector space equipped with a norm,
and that $Y \subseteq V$.  For continuous mappings $f, g : X \to Y$,
consider
\begin{equation}
\label{H(x, t) = (1 - t) f(x) + t g(x), 0 le t le 1}
	H(x, t) = (1 - t) \, f(x) + t \, g(x),	\quad 0 \le t \le 1.
\end{equation}

	This defines a homotopy between $f$ and $g$ as mappings from
$X$ into $V$.  If $Y$ is convex, then $H$ is a homotopy between $f$
and $g$ as mappings from $X$ into $Y$.

	Depending on the circumstances, it may be that $H(x, t) \in Y$
for every $x \in X$ and $0 \le t \le 1$, and hence that $H$ is a
homotopy between $f$ and $g$ as mappings from $X$ into $Y$.

\section{Contractability}
\label{contractability}
\setcounter{equation}{0}

	A topological space $X$ is said to be \emph{contractable} if
there is a continuous mapping $H : X \times [0, 1] \to X$ such that
$H(x, 0) = x$ for every $x \in X$ and $H(x, 1)$ is constant.  More
precisely, this means that there is a $p \in X$ such that $H(x, 1) =
p$ for every $x \in X$.

	Equivalently, $X$ is contractable if the identity mapping on
$X$ is homotopic to a constant as mappings from $X$ to itself.

	If $X$ is contractable, then every element of $X$ can be
connected to a fixed element of $X$ by a continuous path.  Therefore
$X$ is pathwise connected.

	If $V$ is a real vector space equipped with a norm, then $V$
is contractable.  Convex subsets of $V$ are also contractable.

	Let $X$, $Y$ be topological spaces, and let $f$ be a
continuous mapping from $X$ to $Y$.  We say that $f$ is
\emph{homotopically trivial} if $f$ is homotopic to a constant
mapping.

	If $Y$ is pathwise connected, then all constant mappings from
$X$ to $Y$ are homotopic to each other.

	If $Y$ is contractable, then every continuous mapping $f : X
\to Y$ is homotopically trivial.

	Similarly, if $X$ is contractable, then every continuous
mapping $f : X \to Y$ is homotopically trivial.

\section{Mappings into spheres}
\label{mappings into spheres}
\setcounter{equation}{0}

	Let $n$ be a positive integer, and let ${\bf S}^n$ be the
standard unit sphere in ${\bf R}^{n + 1}$, i.e.,
\begin{equation}
	{\bf S}^n = \bigg\{w \in {\bf R}^{n + 1} :
			\sum_{i = 1}^{n + 1} w_i^2 = 1 \bigg\}.
\end{equation}

	For every $p \in {\bf S}^n$, ${\bf S}^n \backslash \{p\}$ is
homeomorphic to ${\bf R}^n$.  In particular, ${\bf S}^n \backslash
\{p\}$ is contractable.

	Let $X$ be a topological space, and let $f : X \to {\bf S}^n$
be a continuous mapping.  If there is a $p \in {\bf S}^n$ such that
$f(x) \ne p$ for every $x \in X$, then $f$ is homotopically trivial.

	It can happen that $f(X) = {\bf S}^n$, and $f$ is
homotopically trivial.  For example, $f$ might be a continuous mapping
from an interval onto the unit circle.

	If $X$ is sufficiently ``small'', then one expects that any
continuous mapping $f : X \to {\bf S}^n$ is homotopically trivial.

	This is not simply a matter of showing that $f(X) \ne {\bf
S}^n$.  It can happen that $f(X) = {\bf S}^n$ even for relative small
spaces $X$ and continuous mappings $f : X \to {\bf S}^n$.

	If $f : X \to {\bf S}^n$ is sufficiently regular, then $f(X)
\ne {\bf S}^n$ may hold automatically.  A basic trick is to show that
the $n$-dimensional volume of $f(X)$ is equal to $0$ in ${\bf S}^n$.
This holds when $f$ is a continuously-differentiable mapping of ${\bf
S}^k$ into ${\bf S}^n$ and $k < n$, for instance.

	Another basic trick is to approximate continuous mappings by
more regular mappings.  If the approximation is homotopy equivalent to
the original mapping and is homotopically trivial, then the original
mapping is homotopically trivial too.

	Using this approach, one can show that every continuous
mapping from ${\bf S}^k$ into ${\bf S}^n$ is homotopically trivial
when $k < n$.  It is well known that the identity mapping on ${\bf
S}^n$ is not homotopically trivial, which is the same as saying that
${\bf S}^n$ is not contractable.

\section{Homotopy equivalence}
\label{homotopy equivalence}
\setcounter{equation}{0}

	Let $X$, $Y$ be topological spaces, and suppose that $f : X
\to Y$ and $g : Y \to X$ are continuous mappings.  If $g \circ f$ is
homotopy equivalent to the identity mapping on $X$, and $f \circ g$ is
homotopy equivalent to the identity mapping on $Y$, then $f$, $g$ are
\emph{homotopy equivalences} from $X$ to $Y$ and $Y$ to $X$, respectively.

	This includes the case when $g \circ f$ is the identity
mapping on $X$ and $f \circ g$ is the identity mapping on $Y$, i.e.,
when $f$ is a homeomorphism from $X$ onto $Y$ and $g$ is the inverse
mapping.

	As a special case, a topological space is contractable if and
only if it is homotopy equivalent to a space with one element.  The
Cartesian product of a topological space $X$ and a contractable space
is homotopy equivalent to $X$.

	Suppose that $X_1$, $X_2$, and $X_3$ are topological spaces
and that
\begin{equation}
	f_1, f_1' : X_1 \to X_2 \quad\hbox{and}\quad f_2, f_2' : X_2 \to X_3
\end{equation}
are continuous mappings.  If $f_1$ is homotopic to $f_1'$ and $f_2$ is
homotopic to $f_2'$, then $f_2 \circ f_1$ is homotopic to $f_2' \circ
f_1'$.  If $X_1$ is homotopy equivalent to $X_2$ and $X_2$ is homotopy
equivalent to $X_3$, then $X_1$ is homotopy equivalent to $X_3$.

\section{Induced mappings}
\label{induced mappings}
\setcounter{equation}{0}

	Let $X$, $Y$ be topological spaces, and let $[X, Y]$ be the
collection of homotopy classes of continuous mappings from $X$ to $Y$.

	Suppose that $X'$ is another topological space and that $\phi
: X' \to X$ is a continuous mapping.  Consider the mapping $\Phi :
\mathcal{C}(X, Y) \to \mathcal{C}(X', Y)$ defined by
\begin{equation}
	\Phi(f) = f \circ \phi.
\end{equation}

	If $f$, $g$ are homotopic mappings from $X$ to $Y$, then $f
\circ \phi$, $g \circ \phi$ are homotopic mappings from $X'$ to $Y$.
Therefore $\Phi$ induces a mapping from $[X, Y]$ to $[X', Y]$.

	If $\widetilde{\phi} : X' \to X$ is a continuous mapping which
is homotopic to $\phi$, then $f \circ \widetilde{\phi}$ is homotopic
to $f \circ \phi$ for every $f \in \mathcal{C}(X, Y)$.  This means
that $\phi$, $\widetilde{\phi}$ induce the same mappings from $[X, Y]$
to $[X', Y]$.

	Let $X''$ be another topological space, let $\phi' : X'' \to
X'$ be a continuous mapping, and let $\Phi' : \mathcal{C}(X', Y) \to
\mathcal{C}(X'', Y)$ be defined by precomposition with $\phi'$, as
before.  Thus $\Phi' \circ \Phi$ is the same as the mapping from
$\mathcal{C}(X, Y)$ to $\mathcal{C}(X'', Y)$ defined by precomposition
with $\phi \circ \phi'$.  The mapping from $[X, Y]$ to $[X'', Y]$
induced by $\phi \circ \phi'$ is the same as the composition of the
mapping from $[X, Y]$ to $[X', Y]$ induced by $\phi$ and the mapping
from $[X', Y]$ to $[X'', Y]$ induced by $\phi'$.

	Similarly, suppose that $Y'$ is another topological space and
$\psi : Y \to Y'$ is continuous.  Consider the mapping $\Psi :
\mathcal{C}(X, Y) \to \mathcal{C}(X, Y')$ defined by $\Psi(f) = \psi
\circ f$, which induces a mapping from $[X, Y]$ to $[X, Y']$.

	If $\widetilde{\psi}$ is another continuous mapping from $Y$
to $Y'$ which is homotopic to $\psi$, then $\widetilde{\psi} \circ f$
is homotopic to $\psi \circ f$ for every continuous mapping $f$ from
$X$ to $Y$.  This implies again that $\psi$, $\widetilde{\psi}$ induce
the same mappings from $[X, Y]$ to $[X, Y']$.

	If $Y''$ is another topological space, $\psi' : Y' \to Y''$ is
another continuous mapping, and $\Psi' : \mathcal{C}(X, Y') \to
\mathcal{C}(X, Y'')$ is defined by postcomposition with $\psi'$, then
$\Psi' \circ \Psi$ is the same as the mapping from $\mathcal{C}(X, Y)$
to $\mathcal{C}(X, Y'')$ defined by postcomposition with $\psi' \circ
\psi$, and the mapping from $[X, Y]$ to $[X, Y'']$ induced by $\psi'
\circ \psi$ is the same as the composition of the mapping from $[X,
Y]$ to $[X, Y']$ induced by $\psi$ and the mapping from $[X, Y']$ to
$[X, Y'']$ induced by $\psi'$.

	In particular, if $\phi : X' \to X$ is a homotopy equivalence,
then $\phi$ induces a one-to-one mapping from $[X, Y]$ onto $[X', Y]$.
If $\psi : Y \to Y'$ is a homotopy equivalence, then $\psi$ induces a
one-to-one mapping from $[X, Y]$ onto $[X, Y']$.

\section{Retracts}
\label{retracts}
\setcounter{equation}{0}

	Let $E$, $U$, and $Y$ be sets, with $E \subseteq U \subseteq
Y$.  A mapping $r : U \to E$ is a \emph{retract} of $U$ onto $E$ if
$r(y) = y$ for every $y \in E$.

	We shall be interested in this when $Y$ is a topological space
and $r$ is a continuous mapping.  In practice, $E$ might be a closed
set or at least relatively closed in $U$, and $U$ might be an open set
or contain $E$ in its interior.

	Let $p$, $q$ be elements of $E$.  If there is a continuous
path in $E$ that connects $p$ to $q$, then the same path connects $p$
to $q$ in $U$.  Conversely, if there is a continuous path in $U$ that
connects $p$ to $q$, then the retract maps that to a path in $E$
connecting $p$ to $q$.

	Similarly, let $X$ be a topological space, and let $f$, $g$ be
continuous mappings from $X$ into $E$.  If $f$, $g$ are homotopic as
mappings from $X$ into $E$, then they are homotopic as mappings from
$X$ into $U$, and the existence of a retract implies that the converse
holds.

	Suppose that $V$ is a real vector space equipped with a norm
$\|v\|$, $U$ is an open set in $V$, $E \subseteq U$, and $r : U \to E$
is a continuous retract.

	We have seen that $U$ has only finitely or countably many
pathwise connected components when $V$ is separable, and hence the
same statement holds for $E$.

	In any case, $E$ is covered by the pathwise connected
components of $U$.  If $E$ is compact, then $E$ is covered by finitely
many of these components, and $E$ has only finitely many pathwise
connected components.  If $E$ is connected, then $E$ is contained in a
single pathwise connected component of $U$, and $E$ is pathwise connected.

	Continuous mappings into $U$ are homotopic when they are
sufficiently close together, using the linear homotopy (\ref{H(x, t) =
(1 - t) f(x) + t g(x), 0 le t le 1}).  Because of the retract, $E$ has
analogous regularity properties.

	As a simple example,
\begin{equation}
	r(v) = \frac{v}{\|v\|}
\end{equation}
is a continuous retract from $U = V \backslash \{0\}$ onto
\begin{equation}
	E = \{v \in V : \|v\| = 1\}.
\end{equation}

\section{Topological groups}
\label{topological groups}
\setcounter{equation}{0}

	Suppose that $G$ is a topological group, which means that $G$
is a group with a topological structure such that the group operations
are continuous.

	If $X$ is a topological space, then the group operation on $G$
can be applied pointwise to functions on $X$ with values in $G$, which
makes the space $\mathcal{C}(X, G)$ of such functions into a group.

	This leads to a group structure on the space $[X, G]$ of
homotopy classes of continuous mappings from $X$ into $G$ as well.

	If $G$ is commutative, then $\mathcal{C}(X, G)$ and $[X, G]$
are commutative too.

	Let $H$ be another topological group, and $\psi : G \to H$ be
a continuous group homomorphism.  The mapping $\Psi : \mathcal{C}(X,
G) \to \mathcal{C}(X, H)$ defined by $\Psi(f) = \psi \circ f$ is a
group homomorphism.  The induced mapping from $[X, G]$ to $[X, H]$ is
also a group homomorphism.

	If $X'$ is another topological space and $\phi : X' \to X$ a
continuous mapping, then the mapping $\Phi : \mathcal{C}(X, G) \to
\mathcal{C}(X', G)$ defined by $\Phi(f) = f \circ \phi$ is a group
homomorphism.  The induced mapping from $[X, G]$ to $[X', G]$ is also
a group homomorphism.

	There are analogous statements for topological semigroups.

\section{The circle group}
\label{circle group}
\setcounter{equation}{0}

	Let ${\bf Z}$ be the integers, a subgroup of ${\bf R}$ with
respect to addition, and let ${\bf T}$ be the quotient group ${\bf R}
/ {\bf Z}$.  There is a natural topology on ${\bf T}$ for which the
quotient mapping ${\bf R} \to {\bf T}$ is a local homeomorphism.  This
makes ${\bf T}$ a compact commutative topological group which is
isomorphic to the unit circle in the complex plane with respect to
multiplication.

	Let $X$ be a topological space, and let $f : X \to {\bf T}$ be
a continuous mapping.  If $\phi : X \to {\bf R}$ is continuous, then
the composition of $\phi$ with the quotient mapping from ${\bf R}$
onto ${\bf T}$ is a continuous mapping from $X$ into ${\bf T}$.

	Let $\mathcal{U} = \{U_\alpha\}_{\alpha \in A}$ be an open
covering of $X$.  Under suitable conditions, for each $\alpha \in A$
there is a continuous mapping $\phi_\alpha : U_\alpha \to {\bf R}$
such that $f$ is equal to the composition of $\phi_\alpha$ with the
quotient mapping from ${\bf R}$ onto ${\bf T}$ on $U_\alpha$.

	For each $\alpha, \beta \in A$ such that $U_\alpha \cap
U_\beta \ne \emptyset$, $\phi_\alpha - \phi_\beta$ takes values in
${\bf Z}$ on $U_\alpha \cap U_\beta$.  If $U_\alpha \cap U_\beta$ is
connected, then $\phi_\alpha - \phi_\beta$ is constant on $U_\alpha
\cap U_\beta$.

\section{The fundamental group}
\label{fundamental group}
\setcounter{equation}{0}

	Let $X$ be a topological space, and let $b$ be a fixed element
of $X$.  Consider continuous mappings from the unit circle ${\bf S}^1$
into $X$ which send $(1, 0)$ to $b$.  Equivalently, the unit circle
can be represented as ${\bf R} / {\bf Z}$, and one can consider
continuous mappings into $X$ which send $0$ to $b$.

	These continuous mappings are also known as \emph{based loops}
in $X$.  In this context, we are concerned with homotopies which
satisfy the same condition.  This leads to an equivalence relation and
a partition of based loops into equivalence classes.

	There is a natural way to compose two based loops in $X$, by
following one and then the other.  The inverse of a loop is defined by
traversing it backwards.  These operations lead to a group structure
on the set of equivalence classes of based loops in $X$, the
fundamental group of $X$ for the basepoint $b$.

	Suppose that $b' \in X$ is another choice of basepoint, and
that $\alpha$ is a continuous path in $X$ from $b$ to $b'$.  A loop
based at $b$ can be transformed into a loop based at $b'$ by
traversing $\alpha$ backwards, then the loop at $b$, and then
$\alpha$.  This transformation sends homotopies of loops based at $b$
to homotopies of loops based at $b'$, and leads to an isomorphism from
the fundamental group of $X$ based to $b$ to the fundamental group
based at $b'$.  If $X$ is pathwise connected, then the fundamental
group of $X$ at any basepoint is isomorphic to the one at any other
basepoint.

\section{Projective spaces}
\label{projective spaces}
\setcounter{equation}{0}

	Let $V$ be a vector space over a field $k$, like the real
numbers.  By definition, the \emph{projective space} ${\bf P}(V)$
associated to $V$ consists of the lines in $V$ through $0$, which is
to say the one-dimensional linear subspaces of $V$.

	Alternatively, let $\sim$ be the relation on $V \backslash
\{0\}$ defined by $v \sim w$ when $v, w \in V$, $v, w \ne 0$, and
there is a $t \in k$ such that $t \ne 0$ and $w = t \, v$.  It is easy
to check that this is an equivalence relation, and one can interpret
${\bf P}(V)$ as the space of equivalence classes determined by $\sim$.

	Thus there is a canonical quotient mapping from $V \backslash
\{0\}$ onto ${\bf P}(V)$ that sends $v \in V$ with $v \ne 0$ to the
element of ${\bf P}(V)$ corresponding to the line in $V$ through $v$.

	If $W$ is a linear subspace of $V$, then every line in $W$ is
also a line in $V$ and therefore ${\bf P}(W) \subseteq {\bf P}(V)$.
The ${\bf P}(W)$'s for $W \subseteq V$ are known as \emph{projective
subspaces} of ${\bf P}(V)$.  In particular, \emph{projective lines} in
${\bf P}(V)$ are projective subspaces ${\bf P}(W)$ corresponding to
$2$-dimensional linear subspaces $W$ of $V$.

	More generally, suppose that $V_1$, $V_2$ are vector spaces
over $k$.  If $T$ is a one-to-one linear mapping from $V_1$ into
$V_2$, then $T$ maps lines in $V_1$ to lines in $V_2$.  Hence $T$
induces a natural mapping $\widehat{T} : {\bf P}(V_1) \to {\bf
P}(V_2)$.  Because $T$ maps linear subspaces of $V_1$ to linear
subspaces of $V_2$, $\widehat{T}$ sends projective subspaces of ${\bf
P}(V_1)$ to projective subspaces of ${\bf P}(V_2)$.

	Observe that $\widehat{T} : {\bf P}(V_1) \to {\bf P}(V_2)$ is
one-to-one, since $T : V_1 \to V_2$ is one-to-one.  If $T$ maps $V_1$
onto $V_2$, then $\widehat{T}$ maps ${\bf P}(V_1)$ onto ${\bf
P}(V_2)$.  The inverse of $\widehat{T}$ is the same as the mapping
from ${\bf P}(V_2)$ to ${\bf P}(V_1)$ induced by the inverse of $T$.
If $V_1$, $V_2$, $V_3$ are vector spaces over $k$ and $T_1 : V_1 \to
V_2$, $T_2, : V_2 \to V_3$ are linear mappings, then the mapping from
${\bf P}(V_1)$ to ${\bf P}(V_3)$ induced by the composition $T_2 \circ
T_1$ is the same as the composition of the mappings induced by $T_1$
and $T_2$.  As an important special case, invertible linear mappings
on a vector space $V$ yield invertible mappings on ${\bf P}(V)$.

	Let $W$ be a linear subspace of $V$, and let $u$ be an element
of $V$ which is not in $W$.  Let $A$ be the affine space $W + u$,
consisting of vectors of the form $w + u$ with $w \in W$.  For every
$a \in A$ there is a unique line in $V$ that contains $a$.  This
defines a mapping from $A$ into ${\bf P}(V)$.  If $V$ is spanned by
$W$ and $u$, then every line in $V$ is either contained in $W$ or
contains an element of $A$.

\section{Real projective spaces}
\label{real projective spaces}
\setcounter{equation}{0}

	Let $V$ be a real vector space, and let ${\bf P}(V)$ be the
corresponding projective space, as in the previous section.  We can
also consider the space ${\bf P}_+(V)$ of rays in $V$.  For each $v
\in V$ with $v \ne 0$, the \emph{ray} in $V$ passing through $v$ is
the set of vectors of the form $r \, v$, where $r$ is a positive real
number.  If $\sim_+$ is the relation on $V \backslash \{0\}$ defined
by $v \sim_+ w$ when $w = t \, v$ for some $t > 0$, then $\sim_+$ is
an equivalence relation on $V \backslash \{0\}$, and ${\bf P}_+(V)$ is
the same as the collection of equivalence classes of $V \backslash
\{0\}$ determined by $\sim_+$.

	There is a canonical quotient mapping from $V \backslash
\{0\}$ onto ${\bf P}_+(V)$, which sends each $v \in V$ with $v \ne 0$
to the ray in $V$ passing through $v$.  There is also a natural
mapping from ${\bf P}_+(V)$ onto ${\bf P}(V)$, which sends each ray in
$V$ to the line that contains it.  The composition of these two
mappings is the same as the usual quotient mapping from $V \backslash
\{0\}$ onto ${\bf P}(V)$.  The mapping from ${\bf P}_+(V)$ onto ${\bf
P}(V)$ is two-to-one, because every line in $V$ contains exactly two
rays in $V$.

	Injective linear mappings between real vector spaces determine
mappings between the corresponding spaces of rays, in the same way as
for projective spaces, and with similar properties.  It is easy to see
that these induced mappings between spaces of rays are compatible with
the induced mappings between the corresponding projective spaces, in
the following sense.  Suppose that $V_1$, $V_2$ be real vector spaces,
and that $T$ is a one-to-one linear mapping from $V_1$ into $V_2$.
The composition of the mapping from ${\bf P}_+(V_1)$ into ${\bf P}_+(V_2)$
induced by $T$ with the natural mapping from ${\bf P}_+(V_2)$ onto
${\bf P}(V_2)$ is the same as the composition of the natural mapping
from ${\bf P}_+(V_1)$ onto ${\bf P}(V_1)$ composed with the mapping
from ${\bf P}(V_1)$ into ${\bf P}(V_2)$ induced by $T$.

	Let $V$ be a real vector space equipped with a norm $\|v\|$,
and let $\Sigma$ be the corresponding unit sphere in $V$, i.e.,
$\Sigma = \{v \in V : \|v\| = 1\}$.  The restriction of the quotient
mapping from $V \backslash \{0\}$ onto ${\bf P}_+(V)$ to $\Sigma$ is a
one-to-one mapping from $\Sigma$ onto ${\bf P}_+(V)$.  The restriction
of the quotient mapping from $V \backslash \{0\}$ onto ${\bf P}(V)$ to
$\Sigma$ is a two-to-one mapping from $\Sigma$ onto ${\bf P}(V)$.

	Let $n$ be a positive integer, and let ${\bf RP}^n$ be the
projective space ${\bf P}({\bf R}^{n + 1})$ associated to ${\bf R}^{n
+ 1}$ as a real vector space.

	For each hyperplane $H$ in ${\bf R}^{n + 1}$ with $0 \not\in
H$, there is an embedding of $H$ into ${\bf RP}^n$ which sends $x \in
H$ to the line in ${\bf R}^{n + 1}$ through $x$.  There is a natural
topology on ${\bf RP}^n$ for which these embeddings are homeomorphisms
onto open subsets of ${\bf RP}^n$.

	This topology can also be characterized by the property that
the quotient mapping from ${\bf R}^{n + 1} \backslash \{0\}$ onto
${\bf RP}^n$ is a continuous open mapping.  The usual two-to-one
mapping from the unit sphere ${\bf S}^n$ in ${\bf R}^{n + 1}$ onto
${\bf RP}^n$ is a local homeomorphism.

	Suppose that $X$ is a topological space, and that $\phi$ is a
continuous mapping from $X$ into ${\bf RP}^n$.  If $\Phi : X \to {\bf
S}^n$ is a continuous mapping, then the composition of $\Phi$ with the
usual mapping from ${\bf S}^n$ onto ${\bf RP}^n$ is such a mapping.
Let $\mathcal{U} = \{U_\alpha\}_{\alpha \in A}$ be an open covering of
$X$.  Under suitable conditions, for each $\alpha \in A$ there is a
continuous mapping $\phi_\alpha : U_\alpha \to {\bf S}^n$ whose
composition with the usual mapping from ${\bf S}^n$ onto ${\bf RP}^n$
is the restriction of $\phi$ to $U_\alpha$.  Of course, $-\phi_\alpha$
has the same feature, since $x$ and $-x$ are sent to the same point in
${\bf RP}^n$ for every $x \in {\bf S}^n$.  Let $\alpha$, $\beta$ be
elements of $A$ such that $U_\alpha \cap U_\beta \ne \emptyset$.  For
every $p \in U_\alpha \cap U_\beta$, $\phi_\beta(p)$ is equal to
$\phi_\alpha(p)$ or to $-\phi_\alpha(p)$.  This is because
$\phi_\alpha(p), \phi_\beta(p) \in {\bf S}^n$ both correspond to the
same point $\phi(p) \in {\bf RP}^n$.  Thus $\phi_\beta = r_{\alpha,
\beta} \, \phi_\alpha$ on $U_\alpha \cap U_\beta$ for some continuous
function $r_{\alpha, \beta}$ on $U_\alpha \cap U_\beta$ with values in
$\{+1, -1\}$.  If $U_\alpha \cap U_\beta$ is connected, then
$r_{\alpha, \beta}$ is constant, and $\phi_\alpha = \phi_\beta$ or
$\phi_\alpha = - \phi_\beta$ on $U_\alpha \cap U_\beta$.

\section{Complex projective spaces}
\label{complex projective spaces}
\setcounter{equation}{0}

	A complex number $z$ can be expressed as $x + y \, i$, where
$x$ and $y$ are real numbers and $i^2 = -1$.  In this case, the
\emph{complex conjugate} of $z$ is denoted $\overline{z}$ and defined
to be $x - y \, i$, and the \emph{modulus} of $z$ is denoted $|z|$ and
defined to be $\sqrt{x^2 + y^2}$.  The field of complex numbers is
denoted ${\bf C}$.  For every $z, w \in {\bf C}$,
\begin{equation}
	\overline{z + w} = \overline{z} + \overline{w}
\end{equation}
and
\begin{equation}
	\overline{z \, w} = \overline{z} \, \overline{w}.
\end{equation}
Hence
\begin{equation}
	|z \, w| = |z| \, |w|,
\end{equation}
since $z \, \overline{z} = |z|^2$.  For each positive integer $n$, the
space ${\bf C}^n$ of $n$-tuples of complex numbers is a vector space
over the complex numbers with respect to coordinatewise addition and
scalar multiplication.  Let ${\bf CP}^n$ be the projective space ${\bf
P}({\bf C}^{n + 1})$ associated to ${\bf C}^{n + 1}$.

	If $H$ is a complex hyperplane in ${\bf C}^{n + 1}$ that does
not contain $0$, then there is an embedding of $H$ into ${\bf CP}^n$
which sends each $v \in H$ to the point in ${\bf CP}^n$ corresponding
to the complex line in ${\bf C}^{n + 1}$ through $v$.  There is a
natural topology on ${\bf CP}^n$ for which these embeddings are
homeomorphisms onto open subsets of ${\bf CP}^n$.  The quotient
mapping from ${\bf C}^{n + 1} \backslash \{0\}$ onto ${\bf CP}^n$ is a
continuous open mapping with respect to this topology.  Using the
aforementioned embedding, one can show that ${\bf CP}^1$ is
homeomorphic to the standard $2$-dimensional sphere ${\bf S}^2$.  This
follows by identifying ${\bf CP}^1$ with the one-point
compactification of the complex plane.

	For $v = (v_1, \ldots, v_{n + 1}) \in {\bf C}^{n + 1}$, the
Euclidean norm of $v$ is defined to be $(|v_1|^2 + \cdots + |v_{n +
1}|^2)^{1/2}$.  The corresponding unit sphere $\Sigma_n$ in ${\bf
C}^{n + 1}$ consists of vectors with norm equal to $1$ and can be
identified with the ordinary unit sphere ${\bf S}^{2 n + 1}$ in ${\bf
R}^{2 n + 2}$.  The restriction of the usual quotient mapping from
${\bf C}^{n + 1} \backslash \{0\}$ onto ${\bf CP}^n$ to $\Sigma_n$ is
a continuous mapping onto ${\bf CP}^n$.  Two vectors $v$, $w$ in
$\Sigma_n$ are sent to the same element of ${\bf CP}^n$ if and only if
there is a $\zeta \in {\bf C}$ with $|\zeta| = 1$ such that $w = \zeta
\, v$.

	Suppose that $X$ is a topological space, and that $\phi$ is a
continuous mapping from $X$ into ${\bf CP}^n$.  If $\Phi$ is a
continuous mapping from $X$ into $\Sigma_n$, then the composition of
$\Phi$ with the mapping from $\Sigma_n$ onto ${\bf CP}^n$ described
above is a continuous mapping from $X$ into ${\bf CP}^n$.  Let
$\mathcal{U} = \{U_\alpha\}_{\alpha \in A}$ be an open covering of
$X$.  Under suitable conditions, for each $\alpha \in A$ there is a
continuous mapping $\phi_\alpha : U_\alpha \to \Sigma_n$ whose
composition with the usual mapping from $\Sigma_n$ onto ${\bf CP}^n$
is the same as the restriction of $\phi$ to $U_\alpha$.  If $\alpha,
\beta \in A$ and $U_\alpha \cap U_\beta \ne \emptyset$, then there is
a continuous mapping $\zeta_{\alpha, \beta} : U_\alpha \cap U_\beta
\to {\bf C}$ such that $|\zeta_{\alpha, \beta}(p)| = 1$ and
\begin{equation}
	\phi_\beta(p) = \zeta_{\alpha, \beta}(p) \, \phi_\alpha(p)
\end{equation}
for every $p \in U_\alpha \cap U_\beta$.  Under suitable conditions,
there is a continuous mapping $r_{\alpha, \beta} : U_\alpha \cap
U_\beta \to {\bf R}$ whose composition with the usual mapping from
${\bf R}$ onto the circle group is equal to $\zeta_{\alpha, \beta}$.
Since $\zeta_{\beta, \alpha} = 1/\zeta_{\alpha, \beta}$, we can also
choose $r_{\alpha, \beta}$ such that $r_{\beta, \alpha} = - r_{\alpha,
\beta}$.  Let $\alpha$, $\beta$, and $\gamma$ be elements of $A$ such
that $U_\alpha \cap U_\beta \cap U_\gamma \ne \emptyset$, and consider
$v_{\alpha, \beta, \gamma} : U_\alpha \cap U_\beta \cap U_\gamma \to
{\bf R}$ defined by
\begin{equation}
	v_{\alpha, \beta, \gamma}
		= r_{\alpha, \beta} + r_{\beta, \gamma} + r_{\gamma, \alpha}.
\end{equation}
Because
\begin{equation}
	\phi_\alpha
		= \zeta_{\gamma, \alpha} \, \phi_\gamma
		= \zeta_{\gamma, \alpha} \, \zeta_{\beta, \gamma} \, \phi_\beta
 		= \zeta_{\gamma, \alpha} \, \zeta_{\beta, \gamma}
			\, \zeta_{\alpha, \beta} \, \phi_\alpha
\end{equation}
on $U_\alpha \cap U_\beta \cap U_\gamma$,
\begin{equation}
 \zeta_{\alpha, \beta} \, \zeta_{\beta, \gamma} \, \zeta_{\gamma, \alpha}
		= 1
\end{equation}
on $U_\alpha \cap U_\beta \cap U_\gamma$, which implies that
$v_{\alpha, \beta, \gamma}$ takes values in the integers.  If
$U_\alpha \cap U_\beta \cap U_\gamma$ is connected, then $v_{\alpha,
\beta, \gamma}$ is constant.

\section{Local differences}
\label{local differences}
\setcounter{equation}{0}

	Let $X$ be a topological space, and let $\widehat{X}$ be the
copy of $X$ in the diagonal of $X \times X$,
\begin{equation}
	\widehat{X} = \{ (x, y) \in X \times X : y = x \}.
\end{equation}
Suppose that $W$ is an open set in $X \times X$ such that $\widehat{X}
\subseteq W$.  We may as well suppose also that $W$ is symmetric in
the sense that $(y, x) \in W$ whenever $(x, y) \in W$, since otherwise
we can replace $W$ with its intersection with the set of $(y, x) \in X
\times X$ such that $(x, y) \in W$.  Let $F(x, y)$ be a continuous
real-valued function on $W$, and consider the problem of finding a
continuous real-valued function $f$ on $X$ such that
\begin{equation}
\label{f(x) - f(y) = F(x, y)}
	f(x) - f(y) = F(x, y)
\end{equation}
when $(x, y) \in W$.

	An obvious necessary condition for the existence of such a
function $f$ is that
\begin{equation}
	F(y, x) = - F(x, y),
\end{equation}
and in particular
\begin{equation}
	F(x, x) = 0.
\end{equation}
Another necessary condition is that
\begin{equation}
\label{F(x, y) + F(y, z) + F(z, x) = 0}
	F(x, y) + F(y, z) + F(z, x) = 0
\end{equation}
when $x, y, z \in X$ and $(x, y), (y, z), (z, x) \in W$.  For if $f$
satisfies (\ref{f(x) - f(y) = F(x, y)}), then this reduces to
\begin{equation}
	(f(x) - f(y)) + (f(y) - f(z)) + (f(z) - f(x)) = 0.
\end{equation}

	For each $p \in X$, put
\begin{equation}
	W(p) = \{x \in X : (x, p) \in W\}
\end{equation}
and $f_p(x) = F(x, p)$ when $x \in W(p)$.  Thus $p \in W(p)$ and
\begin{equation}
	f_p(x) - f_p(p) = F(x, p)
\end{equation}
automatically.  The cocycle condition (\ref{F(x, y) + F(y, z) + F(z,
x) = 0}) implies that
\begin{eqnarray}
	f_p(x) - f_p(y) & = & F(x, p) - F(y, p) = F(x, p) + F(p, y)	\\
		& = & - F(y, x) = F(x, y)	\nonumber
\end{eqnarray}
when $x, y \in W(p)$ and $(x, y) \in W$.  If $p, q \in X$, $x \in W(p)
\cap W(q)$, and $(p, q) \in W$, then
\begin{eqnarray}
	f_p(x) - f_q(x) & = & F(x, p) - F(x, q) = F(q, x) + F(x, p)	\\
		& = & - F(p, q) = F(q, p).	\nonumber
\end{eqnarray}

	If $f : X \to {\bf R}$ satisfies (\ref{f(x) - f(y) = F(x,
y)}), then $f + c$ also satisfies (\ref{f(x) - f(y) = F(x, y)}) for
any real number $c$.  Conversely, if $f, \widetilde{f} : X \to {\bf
R}$ both satisfy (\ref{f(x) - f(y) = F(x, y)}), then $f -
\widetilde{f}$ satisfies (\ref{f(x) - f(y) = F(x, y)}) with $F = 0$.
Thus $f - \widetilde{f}$ is locally constant on $X$, and hence
constant if $X$ is connected.  Similar remarks apply to solutions of
(\ref{f(x) - f(y) = F(x, y)}) for $(x, y)$ in an open set in $X \times
X$ that contains $\widehat{X}$ and is contained in $W$.

\section{Local families}
\label{local families}
\setcounter{equation}{0}

	Let $X$ be a topological space, and let $\mathcal{U} =
\{U_\alpha\}_{\alpha \in A}$ be an open covering of $X$.  Suppose that
for each $\alpha \in A$, $\phi_\alpha$ is a continuous real-valued
function on $U_\alpha$.  Suppose also that $\phi_\alpha - \phi_\beta$
is constant on $U_\alpha \cap U_\beta$ when $U_\alpha \cap U_\beta \ne
\emptyset$.  Consider the question of whether there is a continuous
real-valued function $\phi$ on $X$ such that $\phi - \phi_\alpha$ is
constant on $U_\alpha$ for every $\alpha \in A$.

	If $\phi : X \to {\bf R}$ has this property, then $\phi + c$
also has this property for every real number $c$.  Conversely, if
$\phi, \widetilde{\phi} : X \to {\bf R}$ both have this property, then
the restriction of $\phi - \widetilde{\phi}$ to $U_\alpha$ is constant
for every $\alpha \in A$.  This implies that $\phi - \widetilde{\phi}$
is locally constant on $X$, and therefore constant when $X$ is
connected.

	For $\alpha, \beta \in A$ with $U_\alpha \cap U_\beta \ne
\emptyset$, let $r_{\alpha, \beta} \in {\bf R}$ be the value of
$\phi_\alpha - \phi_\beta$ on $U_\alpha \cap U_\beta$.  Thus
$r_{\beta, \alpha} = - r_{\alpha, \beta}$, and $r_{\alpha, \alpha} =
0$ in particular.  If $\alpha, \beta, \gamma \in A$ and
\begin{equation}
	U_\alpha \cap U_\beta \cap U_\gamma \ne \emptyset,
\end{equation}
then
\begin{equation}
	r_{\alpha, \beta} + r_{\beta, \gamma} + r_{\gamma, \alpha} = 0.
\end{equation}
This is because
\begin{equation}
	(\phi_\alpha - \phi_\beta) + (\phi_\beta - \phi_\gamma)
		+ (\phi_\gamma - \phi_\alpha) = 0
\end{equation}
on $U_\alpha \cap U_\beta \cap U_\gamma$.

	If there is a $\phi : X \to {\bf R}$ such that $\phi -
\phi_\alpha$ is constant on $U_\alpha$ for every $\alpha \in A$, and
$t_\alpha \in {\bf R}$ is this value, then $r_{\alpha, \beta} =
t_\beta - t_\alpha$ when $U_\alpha \cap U_\beta \ne \emptyset$.
Conversely, suppose that there are real numbers $t_\alpha$ for $\alpha
\in A$ which satisfy this equation.  If $\psi_\alpha = \phi_\alpha +
t_\alpha$, then $\psi_\alpha = \psi_\beta$ on $U_\alpha \cap U_\beta$
when this set is nonempty, and $\phi : X \to {\bf R}$ can be defined
by $\phi = \psi_\alpha$ on $U_\alpha$ for each $\alpha$.

\section{Torsion}
\label{torsion}
\setcounter{equation}{0}

	Let $n$ be a positive integer, and let $n {\bf Z}$ be the set
of integer multiples of $n$.  This is a subgroup of the group ${\bf
Z}$ of integers under addition, and the quotient group ${\bf Z} / n
{\bf Z}$ is a finite abelian group of order $n$.

	Let $X$ be a topological space, and let $\mathcal{U} =
\{U_\alpha\}_{\alpha \in A}$ be an open covering of $X$.  Suppose that
for each $\alpha, \beta \in A$ with $U_\alpha \cap U_\beta \ne
\emptyset$, $\rho_{\alpha, \beta}$ is an element of ${\bf Z} / n {\bf
Z}$ such that
\begin{equation}
\label{rho_{beta, alpha} = - rho_{alpha, beta}}
	\rho_{\beta, \alpha} = - \rho_{\alpha, \beta}
\end{equation}
and
\begin{equation}
\label{rho_{alpha, beta} + rho_{beta, gamma} + rho_{gamma, alpha} = 0}
	\rho_{\alpha, \beta} + \rho_{\beta, \gamma} + \rho_{\gamma, \alpha} = 0
\end{equation}
when $U_\alpha \cap U_\beta \cap U_\gamma \ne \emptyset$.

	For instance, it may be that for each $\alpha \in A$ there is
a $\tau_\alpha \in {\bf Z} / n {\bf Z}$ such that $\rho_{\alpha,
\beta} = \tau_\beta - \tau_\alpha$ when $U_\alpha \cap U_\beta \ne
\emptyset$.  Another possibility is that there are $r_{\alpha, \beta}
\in {\bf Z}$ for $\alpha, \beta \in A$ with $U_\alpha \cap U_\beta \ne
\emptyset$ such that
\begin{equation}
\label{r_{beta, alpha} = - r_{alpha, beta}}
	r_{\beta, \alpha} = - r_{\alpha, \beta}
\end{equation}
and
\begin{equation}
\label{r_{alpha, beta} + r_{beta, gamma} + r_{gamma, alpha} = 0}
	r_{\alpha, \beta} + r_{\beta, \gamma} + r_{\gamma, \alpha} = 0
\end{equation}
when $U_\alpha \cap U_\beta \cap U_\gamma \ne \emptyset$, and for
which $\rho_{\alpha, \beta}$ is the projection of $r_{\alpha, \beta}$
in ${\bf Z} / n {\bf Z}$ by the standard quotient mapping.  If there
are $\tau_\alpha \in {\bf Z} / n {\bf Z}$ such that $\rho_{\alpha,
\beta} = \tau_\beta - \tau_\alpha$ when $U_\alpha \cap U_\beta \ne
\emptyset$, then we can choose $t_\alpha \in {\bf Z}$ which project to
$\tau_\alpha$ in ${\bf Z} / n {\bf Z}$ for each $\alpha \in A$ and
take $r_{\alpha, \beta} = t_\beta - t_\alpha$ when $U_\alpha \cap
U_\beta \ne \emptyset$.

	We can also start with integers $r_{\alpha, \beta}$ for
$\alpha, \beta \in A$ with $U_\alpha \cap U_\beta \ne \emptyset$ which
satisfy (\ref{r_{beta, alpha} = - r_{alpha, beta}}) and
(\ref{r_{alpha, beta} + r_{beta, gamma} + r_{gamma, alpha} = 0}), and
consider their projections $\rho_{\alpha, \beta}$ in ${\bf Z} / n {\bf Z}$.
Thus the $\rho_{\alpha, \beta}$'s automatically satisfy
(\ref{rho_{beta, alpha} = - rho_{alpha, beta}}) and (\ref{rho_{alpha,
beta} + rho_{beta, gamma} + rho_{gamma, alpha} = 0}).  Suppose that
there are $\tau_\alpha \in {\bf Z} / n {\bf Z}$ for each $\alpha \in
A$ such that $\rho_{\alpha, \beta} = \tau_\beta - \tau_\alpha$ when
$U_\alpha \cap U_\beta \ne \emptyset$, and choose $t_\alpha \in {\bf
Z}$ for each $\alpha \in A$ whose projection in ${\bf Z} / n {\bf Z}$
is $\tau_\alpha$.  If $r'_{\alpha, \beta} = r_{\alpha, \beta} +
t_\alpha - t_\beta$, then $r'_{\alpha, \beta}$ projects to $0$ in
${\bf Z} / n {\bf Z}$, which means that $r'_{\alpha, \beta} \in n {\bf
Z}$.

	Suppose that $t_\alpha$, $\alpha \in A$, are integers such
that
\begin{equation}
	r_{\alpha, \beta} = t_\beta - t_\alpha \in n {\bf Z}
\end{equation}
when $U_\alpha \cap U_\beta \ne \emptyset$.  Let $\tau_\alpha$ be the
projection of $t_\alpha$ in ${\bf Z} / n {\bf Z}$ for each $\alpha \in
A$.  Because $r_{\alpha, \beta} \in n {\bf Z}$, $\tau_\beta -
\tau_\alpha = 0$ when $U_\alpha \cap U_\beta \ne \emptyset$.  This
implies that the $\tau_\alpha$'s are all the same when $X$ is
connected, and hence there is a $k \in {\bf Z}$ whose projection in
${\bf Z} / n {\bf Z}$ is equal to the common value of the
$\tau_\alpha$'s.  If $t'_\alpha = t_\alpha - k$, then $t'_\alpha \in n
{\bf Z}$ and
\begin{equation}
	r_{\alpha, \beta} = t'_\beta - t'_\alpha
\end{equation}
when $U_\alpha \cap U_\beta \ne \emptyset$.

	Suppose as before that for $\alpha, \beta \in A$ with
$U_\alpha \cap U_\beta \ne \emptyset$ we have $\rho_{\alpha, \beta}
\in {\bf Z} / n {\bf Z}$ satisfying (\ref{rho_{beta, alpha} = -
rho_{alpha, beta}}) and (\ref{rho_{alpha, beta} + rho_{beta, gamma} +
rho_{gamma, alpha} = 0}).  Let $r_{\alpha, \beta}$ be an integer for
every such $\alpha, \beta \in A$ which projects to $\rho_{\alpha,
\beta}$ in ${\bf Z} / n {\bf Z}$ and satisfies $r_{\beta, \alpha} = -
r_{\alpha, \beta}$.  If $\alpha, \beta, \gamma \in A$ and $U_\alpha
\cap U_\beta \cap U_\gamma \ne \emptyset$, then $r_{\alpha, \beta} +
r_{\beta, \gamma} + r_{\gamma, \alpha}$ projects to $0$ in ${\bf Z} /
n {\bf Z}$, by (\ref{rho_{alpha, beta} + rho_{beta, gamma} +
rho_{gamma, alpha} = 0}).  Therefore $v_{\alpha, \beta, \gamma} =
(r_{\alpha, \beta} + r_{\beta, \gamma} + r_{\gamma, \alpha}) / n$ is
an integer.

	Suppose in addition that there are $l_{\alpha, \beta} \in {\bf
Z}$ for $\alpha, \beta \in A$ with $U_\alpha \cap U_\beta \ne
\emptyset$ such that $l_{\beta, \alpha} = - l_{\alpha, \beta}$ and
\begin{equation}
	v_{\alpha, \beta, \gamma}
		= l_{\alpha, \beta} + l_{\beta, \gamma} + l_{\gamma, \alpha}
\end{equation}
when $U_\alpha \cap U_\beta \cap U_\gamma \ne \emptyset$.  If
\begin{equation}
	r'_{\alpha, \beta} = r_{\alpha, \beta} - n \, l_{\alpha, \beta},
\end{equation}
then $r'_{\alpha, \beta}$ satisfies the analogues of (\ref{r_{beta,
alpha} = - r_{alpha, beta}}) and (\ref{r_{alpha, beta} + r_{beta,
gamma} + r_{gamma, alpha} = 0}).  Moreover, the projection of
$r'_{\alpha, \beta}$ in ${\bf Z} / n {\bf Z}$ is equal to
$\rho_{\alpha, \beta}$.

\section{Another variant}
\label{another variant}
\setcounter{equation}{0}

	Let $X$ be a topological space, and suppose that $\mathcal{U}
= \{U_\alpha\}_{\alpha \in A}$ is an open covering of $X$.

	Suppose that $t_\alpha$, $\alpha \in A$, are real numbers such
that $r_{\alpha, \beta} = t_\beta - t_\alpha \in {\bf Z}$ when
$U_\alpha \cap U_\beta \ne \emptyset$.  If $\tau_\alpha$ is the
projection of $t_\alpha$ in ${\bf R} / {\bf Z}$, then $\tau_\alpha =
\tau_\beta$ when $U_\alpha \cap U_\beta \ne \emptyset$.  Hence
$\tau_\alpha$'s are all the same when $X$ is connected, and there is a
$y \in {\bf R}$ whose projection to ${\bf R} / {\bf Z}$ is the common
value of the $\tau_\alpha$'s.  If $t'_\alpha = t_\alpha - y$, then
$t'_\alpha \in {\bf Z}$ and $r_{\alpha, \beta} = t'_\beta - t'_\alpha$
when $U_\alpha \cap U_\beta \ne \emptyset$.

	Let $r_{\alpha, \beta}$ be real numbers for $\alpha, \beta \in
A$ with $U_\alpha \cap U_\beta \ne \emptyset$.  Suppose that there are
$\tau_\alpha \in {\bf R} / {\bf Z}$ for $\alpha \in A$ such that the
projection of $r_{\alpha, \beta}$ in ${\bf R} / {\bf Z}$ is equal to
the difference between $\tau_\beta$ and $\tau_\alpha$ when $U_\alpha
\cap U_\beta \ne \emptyset$.  For each $\alpha \in A$, let $t_\alpha$
be a real number whose projection to ${\bf R} / {\bf Z}$ is
$\tau_\alpha$.  If $r'_{\alpha, \beta} = r_{\alpha, \beta} + t_\alpha
- t_\beta$ when $U_\alpha \cap U_\beta \ne \emptyset$, then
$r'_{\alpha, \beta} \in {\bf Z}$.

	Now let $\rho_{\alpha, \beta}$ be elements of ${\bf R} / {\bf
Z}$ for $\alpha, \beta \in A$ with $U_\alpha \cap U_\beta \ne
\emptyset$, where $\rho_{\beta, \alpha}$ is the inverse of
$\rho_{\alpha, \beta}$ in ${\bf R} / {\bf Z}$, and the combination of
$\rho_{\alpha, \beta}$, $\rho_{\beta, \gamma}$, and $\rho_{\gamma,
\alpha}$ is $0$ in ${\bf R} / {\bf Z}$ when $U_\alpha \cap U_\beta
\cap U_\gamma \ne \emptyset$.  If there are $\tau_\alpha \in {\bf R} /
{\bf Z}$ for each $\alpha \in A$ such that $\rho_{\alpha, \beta}$ is
the difference of $\tau_\beta$ and $\tau_\alpha$ in ${\bf R} / {\bf
Z}$ when $U_\alpha \cap U_\beta \ne \emptyset$, then we can choose
$t_\alpha \in {\bf R}$ for each $\alpha \in A$ such that the
projection of $t_\alpha$ in ${\bf R} / {\bf Z}$ is equal to
$\tau_\alpha$, and $\rho_{\alpha, \beta}$ is the projection of
$t_\beta - t_\alpha$ in ${\bf R} / {\bf Z}$.  Otherwise, choose
$r_{\alpha, \beta} \in {\bf R}$ for $\alpha, \beta \in A$ with
$U_\alpha \cap U_\beta \ne \emptyset$ such that $r_{\beta, \alpha} = -
r_{\alpha, \beta}$ and $\rho_{\alpha, \beta}$ is the projection of
$r_{\alpha, \beta}$ in ${\bf R} / {\bf Z}$.  Put
\begin{equation}
	v_{\alpha, \beta, \gamma}
		= r_{\alpha, \beta} + r_{\beta, \gamma} + r_{\gamma, \alpha}
\end{equation}
when $U_\alpha \cap U_\beta \cap U_\gamma \ne \emptyset$.  The
projection of $v_{\alpha, \beta, \gamma}$ in ${\bf R} / {\bf Z}$ is
$0$, and therefore $v_{\alpha, \beta, \gamma} \in {\bf Z}$.  Suppose
that there are $l_{\alpha, \beta} \in {\bf Z}$ for $\alpha, \beta \in
A$ with $U_\alpha \cap U_\beta \ne \emptyset$ such that $l_{\beta,
\alpha} = - l_{\alpha, \beta}$ and
\begin{equation}
	v_{\alpha, \beta, \gamma}
		= l_{\alpha, \beta} + l_{\beta, \gamma} + l_{\gamma, \alpha}
\end{equation}
when $U_\alpha \cap U_\beta \cap U_\gamma \ne \emptyset$.  If
\begin{equation}
	r'_{\alpha, \beta} = r_{\alpha, \beta} - l_{\alpha, \beta},
\end{equation}
then $r'_{\beta, \alpha} = - r'_{\alpha, \beta}$,
\begin{equation}
	r'_{\alpha, \beta} + r'_{\beta, \gamma} + r'_{\gamma, \alpha} = 0
\end{equation}
when $U_\alpha \cap U_\beta \cap U_\gamma \ne \emptyset$, and
$\rho_{\alpha, \beta}$ is the projection of $r'_{\alpha, \beta}$ in
${\bf R} / {\bf Z}$.

\section{Functions}
\label{functions}
\setcounter{equation}{0}

	Let $X$ be a topological space, and let $\mathcal{U} =
\{U_\alpha\}_{\alpha \in A}$ be an open covering of $X$.  Suppose that
$\mathcal{U}$ is locally finite on $X$, in the sense that every point
in $X$ has a neighborhood which intersects $U_\alpha$ for only
finitely many $\alpha \in A$.  Under suitable conditions, there is a
partition of unity subordinate to $\mathcal{U}$ on $X$, which is to
say a collection of continuous real-valued functions $\eta_\alpha$ on
$X$ for $\alpha \in A$ such that $0 \le \eta_\alpha(p) \le 1$ for
every $\alpha \in A$ and $p \in X$, $\eta_\alpha(q) = 0$ for every $q$
in a neighborhood of any $p \in X \backslash U_\alpha$, and
\begin{equation}
	\sum_{\alpha \in A} \eta_\alpha(p) = 1
\end{equation}
for every $p \in X$.  Note that there are only finitely many nonzero
terms in the sum for each $p \in X$.

	Suppose that for each $\alpha, \beta \in A$ with $U_\alpha
\cap U_\beta \ne \emptyset$ there is a continuous real-valued function
$r_{\alpha, \beta}$ on $U_\alpha \cap U_\beta$ such that $r_{\beta,
\alpha} = - r_{\alpha, \beta}$ and
\begin{equation}
	r_{\alpha, \beta} + r_{\beta, \gamma} + r_{\gamma, \alpha} = 0
\end{equation}
on $U_\alpha \cap U_\beta \cap U_\gamma$ when this set is not empty.

	For every $\alpha \in A$ and $p \in U_\alpha$, let
$t_\alpha(p)$ be the sum of $- r_{\alpha, \gamma}(p) \, \eta_\gamma(p)$
over $\gamma \in A$ such that $p \in U_\gamma$.  Basically,
\begin{equation}
	t_\alpha(p) = - \sum_\gamma r_{\alpha, \gamma}(p) \, \eta_\gamma(p),
\end{equation}
where $r_{\alpha, \gamma}(p) \, \eta_\gamma(p)$ is interpreted as
being $0$ when $p \not\in U_\gamma$, even though $r_{\alpha,
\gamma}(p)$ is not defined.  This is a continuous real-valued function
on $U_\alpha$.

	If $p \in U_\alpha \cap U_\beta$, then
\begin{eqnarray}
	t_\beta(p) - t_\alpha(p) & = & 
 \sum_\gamma (r_{\alpha, \gamma}(p) - r_{\beta, \gamma}(p)) \, \eta_\gamma(p)
								\\
	& = & \sum_\gamma r_{\alpha, \beta}(p) \, \eta_\gamma(p)
		= r_{\alpha, \beta}(p),				\nonumber
\end{eqnarray}
since only $p \in U_\gamma$ are important for the sum.

\section{More functions}
\label{more functions}
\setcounter{equation}{0}

	Let $X$ be a topological space, and let $\mathcal{U} =
\{U_\alpha\}_{\alpha \in A}$ be an open covering of $X$.  Suppose that
for $\alpha, \beta \in A$ with $U_\alpha \cap U_\beta \ne \emptyset$,
$\rho_{\alpha, \beta}$ is a continuous mapping from $U_\alpha \cap
U_\beta$ into the circle group ${\bf T} = {\bf R} / {\bf Z}$ such that
$\rho_{\beta, \alpha}$ is the inverse of $\rho_{\alpha, \beta}$ in
${\bf T}$ and the combination of $\rho_{\alpha, \beta}$, $\rho_{\beta,
\gamma}$, and $\rho_{\gamma, \alpha}$ is $0$ on $U_\alpha \cap U_\beta
\cap U_\gamma$ when this set is nonempty.  Under suitable conditions,
there are continuous functions $r_{\alpha, \beta} : U_\alpha \cap
U_\beta \to {\bf R}$ such that $r_{\beta, \alpha} = - r_{\alpha,
\beta}$ and $\rho_{\alpha, \beta}$ is the composition of $r_{\alpha,
\beta}$ with the usual quotient mapping from ${\bf R}$ onto ${\bf T}$.
Put
\begin{equation}
	v_{\alpha, \beta, \gamma}
		= r_{\alpha, \beta} + r_{\beta, \gamma} + r_{\gamma, \alpha}
\end{equation}
on $U_\alpha \cap U_\beta \cap U_\gamma$ when this set is nonempty.
By hypothesis, $v_{\alpha, \beta, \gamma}$ takes values in the
integers, and hence is constant if $U_\alpha \cap U_\beta \cap
U_\gamma$ is connected.  Suppose that there are $l_{\alpha, \beta} \in
{\bf Z}$ for $\alpha, \beta \in A$ with $U_\alpha \cap U_\beta \ne
\emptyset$ such that $l_{\beta, \alpha} = - l_{\alpha, \beta}$ and
\begin{equation}
	v_{\alpha, \beta, \gamma}
		= l_{\alpha, \beta} + l_{\beta, \gamma} + l_{\gamma, \alpha}
\end{equation}
when $U_\alpha \cap U_\beta \cap U_\gamma \ne \emptyset$.
If $r'_{\alpha, \beta} = r_{\alpha, \beta} - l_{\alpha, \beta}$, then
$\rho_{\alpha, \beta}$ is the composition of $r'_{\alpha, \beta}$ with
the quotient mapping from ${\bf R}$ onto ${\bf T}$, and
\begin{equation}
	r'_{\alpha, \beta} + r'_{\beta, \gamma} + r'_{\gamma, \alpha} = 0
\end{equation}
on $U_\alpha \cap U_\beta \cap U_\gamma$ when this set is nonempty.

\section{Circle bundles}
\label{circle bundles}
\setcounter{equation}{0}

	Let $X$ be a topological space, and let $\mathcal{U} =
\{U_\alpha\}_{\alpha \in A}$ be an open covering of $X$.  Suppose that
for every $\alpha, \beta \in W$ with $U_\alpha \cap U_\beta \ne
\emptyset$, $\rho_{\alpha, \beta}$ is a continuous mapping from
$U_\alpha \cap U_\beta$ to the circle group ${\bf T} = {\bf R} / {\bf
Z}$ such that $\rho_{\beta, \alpha}$ is the inverse of $\rho_{\alpha,
\beta}$ and the combination of $\rho_{\alpha, \beta}$, $\rho_{\beta,
\gamma}$, and $\rho_{\gamma, \alpha}$ is $0$ on $U_\alpha \cap U_\beta
\cap U_\gamma$ when this set is nonempty.

	For each $\alpha \in A$, let $E_\alpha$ be the Cartesian
product of $U_\alpha$ and ${\bf T}$.  Let $E$ be the space obtained by
gluing the $E_\alpha$'s together by identifying $(p, z) \in E_\alpha$
with $(p, z')$ when $p \in U_\alpha \cap U_\beta$ and $z'$ is the
combination of $z$ and $\rho_{\alpha, \beta}(p)$ in ${\bf T}$.  The
conditions on the $\rho_{\alpha, \beta}$'s imply that these
identifications are consistent with each other.  The obvious
projections from $E_\alpha$ to $U_\alpha$ lead to a continuous mapping
from $E$ onto $X$.  For each $p \in X$, this mapping sends a copy of
${\bf T}$ in $E$ onto $p$.

	Suppose that $\Phi : X \to E$ is a continuous mapping whose
composition with the natural projection from $E$ onto $X$ is the
identity mapping.  For each $\alpha \in A$, the restriction of $\Phi$
to $U_\alpha$ corresponds to a mapping $\Phi_\alpha : U_\alpha \to
E_\alpha$ of the form $\Phi_\alpha(p) = (p, \phi_\alpha(p))$ for some
continuous mapping $\phi_\alpha : U_\alpha \to {\bf T}$.  Because of
the identifications in $E$, $\phi_\beta(p)$ ought to be the
combination of $\rho_{\alpha, \beta}(p)$ and $\phi_\alpha(p)$ in ${\bf
T}$ when $p \in U_\alpha \cap U_\beta$.  Conversely, if for each
$\alpha \in A$ there is a continuous mapping $\phi_\alpha : U_\alpha
\to {\bf T}$ such that $\phi_\beta$ is the combination of
$\rho_{\alpha, \beta}$ and $\phi_\alpha$ on $U_\alpha \cap U_\beta$
when this set is nonempty, then one can reverse the process and get a
continuous mapping $\Phi : X \to E$ whose composition with the
projection from $E$ onto $X$ is the identity mapping.  One can also
show that $E$ is equivalent to $X \times {\bf T}$.

	It is sometimes convenient to think of ${\bf T}$ as the unit
circle in the complex plane with multiplication as the group
operation.  One can use the $\rho_{\alpha, \beta}$'s to glue $U_\alpha
\times {\bf C}$ for $\alpha \in A$ in the same way to get a complex
line bundle $L$ over $X$.  There is automatically a continuous mapping
from $X$ into $L$ which projects to the indentity mapping on $X$,
using the element that corresponds to $0$ in each copy of ${\bf C}$.
Continuous mappings from $X$ to nonzero elements of $L$ that project
to the identity mapping on $X$ are like maps from $X$ to $E$ in the
previous paragraph.

\section{Noncommutative groups}
\label{noncommutative groups}
\setcounter{equation}{0}

	Let $G$ be a group, which may not be commutative.  We may also
ask $G$ to be a topological group.  Let $X$ be a topological space,
and let $\mathcal{U} = \{U_\alpha\}_{\alpha \in A}$ be an open
covering of $X$.

	Suppose that for every $\alpha, \beta \in A$ with $U_\alpha
\cap U_\beta \ne \emptyset$ we have $\rho_{\alpha, \beta}$ which is an
element of $G$ or a continuous mapping from $U_\alpha \cap U_\beta$
into $G$ when $G$ is a topological group such that
\begin{equation}
\label{rho_{beta, alpha} = rho_{beta, alpha}^{-1}}
	\rho_{\beta, \alpha} = \rho_{\beta, \alpha}^{-1}
\end{equation}
and
\begin{equation}
\label{rho_{alpha, beta} rho_{beta, gamma} rho_{gamma, alpha} = e}
 \rho_{\alpha, \beta} \, \rho_{\beta, \gamma} \, \rho_{\gamma, \alpha} = e
\end{equation}
when $U_\alpha \cap U_\beta \cap U_\gamma \ne \emptyset$, where $e$ is
the identity element of $G$.  Alternatively, one might have
\begin{equation}
\label{rho_{gamma, alpha} rho_{beta, gamma} rho_{alpha, beta} = e}
 \rho_{\gamma, \alpha} \, \rho_{\beta, \gamma} \, \rho_{\alpha, \beta} = e
\end{equation}
in place of (\ref{rho_{alpha, beta} rho_{beta, gamma} rho_{gamma,
alpha} = e}).

	Of course, any group is a topological group with respect to
the discrete topology.  Any element of a group can be considered as a
constant mapping into the group.  If the $\rho_{\alpha, \beta}$'s are
functions, then (\ref{rho_{beta, alpha} = rho_{beta, alpha}^{-1}}) is
supposed to hold at every point in $U_\alpha \cap U_\beta$, and
(\ref{rho_{alpha, beta} rho_{beta, gamma} rho_{gamma, alpha} = e}) or
(\ref{rho_{gamma, alpha} rho_{beta, gamma} rho_{alpha, beta} = e}) is
supposed to hold at every point in $U_\alpha \cap U_\beta \cap
U_\gamma$.

	If $\tau_\alpha$ is an element of $G$ or a continuous mapping
from $U_\alpha$ into $G$ when $G$ is a topological group, then
$\rho_{\alpha, \beta} = \tau_\alpha^{-1} \, \tau_\beta$ automatically
satisfies (\ref{rho_{beta, alpha} = rho_{beta, alpha}^{-1}}) and
(\ref{rho_{alpha, beta} rho_{beta, gamma} rho_{gamma, alpha} = e}).
Similarly, $\rho_{\alpha, \beta} = \tau_\beta \, \tau_\alpha^{-1}$
satisfies (\ref{rho_{beta, alpha} = rho_{beta, alpha}^{-1}}) and
(\ref{rho_{gamma, alpha} rho_{beta, gamma} rho_{alpha, beta} = e}).
In general, $\rho_{\alpha, \beta}$ satisfies (\ref{rho_{beta, alpha} =
rho_{beta, alpha}^{-1}}) and (\ref{rho_{alpha, beta} rho_{beta, gamma}
rho_{gamma, alpha} = e}) if and only if $\rho_{\alpha, \beta}^{-1}$
satisfies (\ref{rho_{beta, alpha} = rho_{beta, alpha}^{-1}}) and
(\ref{rho_{gamma, alpha} rho_{beta, gamma} rho_{alpha, beta} = e}).

\section{Quaternions}
\label{quaternions}
\setcounter{equation}{0}

	A \emph{quaternion} $x$ can be expressed as
\begin{equation}
\label{x = x_1 + x_2 i + x_3 j + x_4 k}
	x = x_1 + x_2 \, i + x_3 \, j + x_4 \, k,
\end{equation}
where
\begin{equation}
	i^2 = j^2 = k^2 = -1
\end{equation}
and
\begin{equation}
	i \, j = - j \, i = k.
\end{equation}
More precisely, the algebra of quaternions ${\bf H}$ is an associative
algebra over the real numbers with multiplicative identity element
$1$.  As a real vector space, ${\bf H}$ has dimension $4$ and $1$,
$i$, $j$, and $k$ form a basis of ${\bf H}$.  Note that
\begin{equation}
	i \, k = - k \, i = - j
\end{equation}
and
\begin{equation}
	j \, k = - k \, j = i.
\end{equation}
Thus ${\bf H}$ is a noncommutative algebra, since $i$, $j$, and $k$
anticommute.

	Let $x \in {\bf H}$ be given as in (\ref{x = x_1 + x_2 i + x_3
j + x_4 k}).  The \emph{conjugate} $x^*$ of $x$ is defined by
\begin{equation}
	x^* = x_1 - x_2 \, i - x_3 \, j - x_4 \, k.
\end{equation}
Clearly $x \mapsto x^*$ is a real-linear mapping on ${\bf H}$, and it
is easy to see that
\begin{equation}
	(x \, y)^* = y^* \, x^*
\end{equation}
for every $x, y \in {\bf H}$.

	The \emph{modulus} $|x|$ of $x$ is defined by
\begin{equation}
	|x| = \sqrt{x_1^2 + x_2^2 + x_3^2 + x_4^2}.
\end{equation}
One can check that
\begin{equation}
	x \, x^* = x^* \, x = |x|^2.
\end{equation}
Hence
\begin{equation}
	|x \, y| = |x| \, |y|
\end{equation}
for every $x, y \in {\bf H}$.  If $x \ne 0$, then
\begin{equation}
	x^{-1} = \frac{x^*}{|x|^2}
\end{equation}
is the multiplicative inverse of $x$.

	The nonzero quaternions form a noncommutative group under
multiplication.  The modulus defines a homomorphism from this group
onto the positive real numbers with respect to multiplication.  The
quaternions of modulus $1$ form a subgroup of the nonzero quaternions
which is the kernel of this homomorphism.

\section{Quaternionic projective spaces}
\label{quaternionic projective spaces}
\setcounter{equation}{0}

	Let $n$ be a positive integer, and let ${\bf H}^{n + 1}$ be
the space of $(n + 1)$-tuples of quaternions.  Thus ${\bf H}^{n + 1}$
is a real vector space with respect to coordinatewise addition and
multiplication by real numbers, and it is also a kind of vector space
over ${\bf H}$ with respect to coordinatewise multiplication on the
left by elements of ${\bf H}$, i.e.,
\begin{equation}
	t \, v = (t \, v_1, \ldots, t \, v_{n + 1})
\end{equation}
when $t \in {\bf H}$ and $v = (v_1, \ldots, v_{n + 1}) \in {\bf H}^{n
+ 1}$.  If $v \ne 0$, then the set of $t \, v$, $t \in {\bf H}$, is a
quaternionic line in ${\bf H}^{n + 1}$.  By definition, the
quaternionic projective space ${\bf HP}^n$ consists of the
quaternionic lines in ${\bf H}^{n + 1}$.

	Equivalently, let $\sim$ be the relation on ${\bf H}^{n + 1}
\backslash \{0\}$ defined by $v \sim w$ when $w = t \, v$ for some $t
\in {\bf H}$.  This defines an equivalence relation on ${\bf H}^{n +
1} \backslash \{0\}$.  One can identify ${\bf HP}^n$ with the set of
equivalence classes determined by $\sim$ in ${\bf H}^{n + 1}$.

	For $v = (v_1, \ldots, v_{n + 1}) \in {\bf H}^{n + 1}$, put
\begin{equation}
	\|v\| = \sqrt{|v_1|^2 + \cdots + |v_{n + 1}|^2},
\end{equation}
where $|v_1|, \ldots, |v_{n + 1}|$ are the moduli of the quaternions
$v_1, \ldots, v_{n + 1}$.  Consider the set
\begin{equation}
	\{v \in {\bf H}^{n + 1} : \|v\| = 1\},
\end{equation}
which can be identified with the standard sphere ${\bf S}^{4 n + 3}$,
using the obvious correspondence between ${\bf H}^{n + 1}$ and ${\bf
R}^{4 n + 4}$.  There is a natural projection from ${\bf H}^{n + 1}
\backslash \{0\}$ onto this sphere, defined by
\begin{equation}
	v \mapsto \frac{v}{\|v\|}.
\end{equation}
This intersection of this sphere with any quaternionic line in ${\bf
H}^{n + 1}$ can be identified with ${\bf S}^3$.

	There is a natural topology on ${\bf HP}^n$ such that the
quotient mapping from ${\bf H}^{n + 1} \backslash \{0\}$ onto ${\bf
HP}^n$ is a continuous open mapping.  The restriction of the quotient
mapping to the unit sphere in ${\bf H}^{n + 1}$ as in the previous
paragraph as also a continuous open mapping.  This topology can also
be described in terms of homeomorphic embeddings of ${\bf H}^n$ onto
open subsets of ${\bf HP}^n$.  In particular, ${\bf HP}^1$ is a
one-point compactification of ${\bf H} \cong {\bf R}^4$, and thus
homeomorphic to ${\bf S}^4$.

\section{Liftings}
\label{liftings}
\setcounter{equation}{0}

	Let $X$ be a topological space, and let $\phi$ be a continuous
mapping from $X$ into ${\bf HP}^n$.  If $\Phi$ is a continuous mapping
from $X$ into ${\bf H}^{n + 1} \backslash \{0\}$, then the composition
of $\Phi$ with the usual quotient mapping from ${\bf H}^{n + 1}$ onto
${\bf HP}^n$ is such a mapping.  Let $\mathcal{U} =
\{U_\alpha\}_{\alpha \in A}$ be an open covering of $X$.  Under
suitable conditions, for each $\alpha \in A$ there is a continuous
mapping $\phi_\alpha$ from $U_\alpha$ into ${\bf H}^{n + 1} \backslash
\{0\}$ whose composition with the quotient mapping from ${\bf H}^{n +
1} \backslash \{0\}$ onto ${\bf HP}^n$ is equal to $\phi$ on
$U_\alpha$.  If $U_\alpha \cap U_\beta \ne \emptyset$, then there is a
continuous mapping $\rho_{\alpha, \beta}$ from $U_\alpha \cap U_\beta$
to ${\bf H} \backslash \{0\}$ such that
\begin{equation}
	\phi_\beta = \rho_{\alpha, \beta} \, \phi_\alpha
\end{equation}
on $U_\alpha \cap U_\beta$.  Clearly $\rho_{\beta, \alpha} =
\rho_{\alpha, \beta}^{-1}$.  If $U_\alpha \cap U_\beta \cap U_\gamma
\ne \emptyset$, then
\begin{equation}
	\phi_\alpha = \rho_{\gamma, \alpha} \, \phi_\gamma
		= \rho_{\gamma, \alpha} \, \rho_{\beta, \gamma} \, \phi_\beta
 		= \rho_{\gamma, \alpha} \, \rho_{\beta, \gamma}
			\, \rho_{\alpha, \beta} \, \phi_\alpha
\end{equation}
and hence
\begin{equation}
 \rho_{\gamma, \alpha} \, \rho_{\beta, \gamma} \, \rho_{\alpha, \beta} = 1
\end{equation}
on $U_\alpha \cap U_\beta \cap U_\gamma$.  Suppose that for each
$\alpha \in A$ there is a continuous mapping $\tau_\alpha$ on
$U_\alpha$ with values in ${\bf H} \backslash \{0\}$ such that
$\rho_{\alpha, \beta} = \tau_\beta \, \tau_\alpha^{-1}$ on $U_\alpha
\cap U_\beta$.  If
\begin{equation}
	\phi'_\alpha = \tau_\alpha^{-1} \, \phi_\alpha,
\end{equation}
then the composition of $\phi'_\alpha$ with the quotient mapping from
${\bf H}^{n + 1} \backslash \{0\}$ onto ${\bf HP}^n$ is also equal to
$\phi$ on $U_\alpha$.  Moreover, $\phi'_\alpha = \phi'_\beta$ on
$U_\alpha \cap U_\beta$ when this set is nonempty.  This permits us to
put $\Phi = \phi'_\alpha$ on $U_\alpha$ for each $\alpha \in X$, to
get a continuous mapping from $X$ into ${\bf H}^{n + 1} \backslash
\{0\}$ whose composition with the quotient mapping is equal to $\phi$.

\section{Homomorphisms}
\label{homomorphisms}
\setcounter{equation}{0}

	Let $G$, $H$ be groups, and let $\phi$ be a homomorphism from
$G$ into $H$.  We may also ask that $G$, $H$ be topological groups,
and that $\phi$ be continuous.  Let $N$ be the normal subgroup of $G$
which is the kernel of $\phi$.  Let $X$ be a topological space, and
let $\mathcal{U} = \{U_\alpha\}_{\alpha \in A}$ be an open covering of
$X$.  Suppose that for $\alpha, \beta \in A$ with $U_\alpha \cap
U_\beta \ne \emptyset$ we have $r_{\alpha, \beta}$ which is an element
of $G$ or a continuous mapping from $U_\alpha \cap U_\beta$ into $G$
when $G$ is a topological group and which satisfies $r_{\beta, \alpha}
= r_{\alpha, \beta}^{-1}$ and
\begin{equation}
	r_{\alpha, \beta} \, r_{\beta, \gamma} \, r_{\gamma, \alpha} = e
\end{equation}
when $U_\alpha \cap U_\beta \cap U_\gamma \ne \emptyset$.  If
$\rho_{\alpha, \beta} = \phi(r_{\alpha, \beta})$, so that
$\rho_{\alpha, \beta}$ is an element of $H$ or a continuous mapping
from $U_\alpha \cap U_\beta$ into $H$, then $\rho_{\beta, \alpha} =
\rho_{\alpha, \beta}^{-1}$ and
\begin{equation}
 \rho_{\alpha, \beta} \, \rho_{\beta, \gamma} \, \rho_{\gamma, \alpha} = e'
\end{equation}
when $U_\alpha \cap U_\beta \cap U_\gamma \ne \emptyset$, where $e'$
denotes the identity element in $H$.  Similar remarks apply to the
cocycle condition in the other order.

	Suppose also that for each $\alpha \in A$ there is a
$\tau_\alpha$ which is an element of $H$ or a continuous mapping from
$U_\alpha$ into $H$ when $H$ is a topological group such that
$\rho_{\alpha, \beta} = \tau_\alpha^{-1} \, \tau_\beta$.  Under
suitable conditions, for each $\alpha \in A$ there is a $t_\alpha$
which is an element of $G$ or a continuous mapping from $U_\alpha$
into $G$ such that $\tau_\alpha = \phi \circ t_\alpha$.  If
$\widetilde{r}_{\alpha, \beta} = t_\alpha \, r_{\alpha, \beta} \,
t_\beta^{-1}$, then $\phi(\widetilde{r}_{\alpha, \beta}) = e'$,
$\widetilde{r}_{\beta, \alpha} = \widetilde{r}_{\alpha, \beta}^{-1}$,
and
\begin{equation}
	\widetilde{r}_{\alpha, \beta} \, \widetilde{r}_{\beta, \gamma}
		\, \widetilde{r}_{\gamma, \alpha} = e
\end{equation}
when $U_\alpha \cap U_\beta \cap U_\gamma \ne \emptyset$.  Thus
$\widetilde{r}_{\alpha, \beta}$ is an element of $N$ or a continuous
mapping from $U_\alpha \cap U_\beta$ into $N$, since
$\phi(\widetilde{r}_{\alpha, \beta}) = e'$.

	Now suppose that we start with $\rho_{\alpha, \beta}$ which is
an element of $H$ or a continuous mapping from $U_\alpha \cap U_\beta$
into $H$ when this set is nonempty such that $\rho_{\beta, \alpha} =
\rho_{\alpha, \beta}^{-1}$ and $\rho_{\alpha, \beta} \, \rho_{\beta,
\gamma} \, \rho_{\gamma, \alpha} = e'$ when $U_\alpha \cap U_\beta
\cap U_\gamma \ne \emptyset$.  If $\tau_\alpha$ and $t_\alpha$ are as
in the previous paragraph, then $r_{\alpha, \beta} = t_\alpha^{-1} \,
t_\beta$ satisfies $\phi(r_{\alpha, \beta}) = \rho_{\alpha, \beta}$,
$r_{\beta, \alpha} = r_{\alpha, \beta}^{-1}$, and $r_{\alpha, \beta}
\, r_{\beta, \gamma} \, r_{\gamma, \alpha} = e$ when $U_\alpha \cap
U_\beta \cap U_\gamma \ne \emptyset$.  Otherwise, under suitable
conditions there are $r_{\alpha, \beta}$ which are elements of $G$ or
continuous mappings from $U_\alpha \cap U_\beta$ into $G$ when this
set is nonempty such that $\phi(r_{\alpha, \beta}) = \rho_{\alpha,
\beta}$ and $r_{\beta, \alpha} = r_{\alpha, \beta}^{-1}$.  If
$v_{\alpha, \beta, \gamma} = r_{\alpha, \beta} \, r_{\beta, \gamma} \,
r_{\gamma, \alpha}$ when $U_\alpha \cap U_\beta \cap U_\gamma \ne
\emptyset$, then $\phi(v_{\alpha, \beta, \gamma}) = e'$, and hence
$v_{\alpha, \beta, \gamma} \in N$.

	Let us make the additional hypothesis that $N$ is contained in
the center of $G$, which is to say that each element of $N$ commutes
with every element of $G$.  Suppose that for each $\alpha, \beta \in
A$ with $U_\alpha \cap U_\beta \ne \emptyset$ there is a
$\lambda_{\alpha, \beta}$ which is an element of $N$ or a continuous
mapping from $U_\alpha \cap U_\beta$ into $N$ such that
$\lambda_{\beta, \alpha} = \lambda_{\alpha, \beta}^{-1}$ and
$v_{\alpha, \beta, \gamma} = \lambda_{\alpha, \beta} \,
\lambda_{\beta, \gamma} \, \lambda_{\gamma, \alpha}$ when $U_\alpha
\cap U_\beta \cap U_\gamma \ne \emptyset$.  If $\widehat{r}_{\alpha,
\beta} = r_{\alpha, \beta} \, \lambda_{\alpha, \beta}^{-1}$, then
$\phi(\widehat{r}_{\alpha, \beta}) = \phi(r_{\alpha, \beta}) =
\rho_{\alpha, \beta}$, $\widehat{r}_{\beta, \alpha} =
\widehat{r}_{\alpha, \beta}^{-1}$, and $\widehat{r}_{\alpha, \beta} \,
\widehat{r}_{\beta, \gamma} \, \widehat{r}_{\gamma, \alpha} = e$ when
$U_\alpha \cap U_\beta \cap U_\gamma \ne \emptyset$.

\section{Vector bundles}
\label{vector bundles}
\setcounter{equation}{0}

	Let $n$ be a positive integer, and let $GL(n, {\bf R})$ be the
group of $n \times n$ invertible matrices with entries in the real
numbers, using matrix multiplication as the group operation.  The $n
\times n$ identity matrix $I$ is the identity element of this group.
Equivalently, $GL(n, {\bf R})$ is the open set of $n \times n$ real
matrices with nonzero determinant, and is a topological group with
respect to the induced topology.  The elements of $GL(n, {\bf R})$
also correspond to invertible linear transformations on ${\bf R}^n$ in
the usual way.

	Let $X$ be a topological space, and let $\mathcal{U} =
\{U_\alpha\}_{\alpha \in A}$ be an open covering of $X$.  Suppose that
for each $\alpha, \beta \in A$ with $U_\alpha \cap U_\beta \ne
\emptyset$, $\rho_{\alpha, \beta}$ is a continuous mapping from
$U_\alpha \cap U_\beta$ into $GL(n, {\bf R})$ such that $\rho_{\beta,
\alpha} = \rho_{\alpha, \beta}$ and
\begin{equation}
 \rho_{\gamma, \alpha} \, \rho_{\beta, \gamma} \, \rho_{\alpha, \beta} = I
\end{equation}
on $U_\alpha \cap U_\beta \cap U_\gamma$ when this set is nonempty.
For each $\alpha \in A$, let $E_\alpha$ be the Cartesian product of
$U_\alpha$ with ${\bf R}^n$.  Let us identify $(p, v) \in E_\alpha$
with $(p, v') \in E_\beta$ when $p \in U_\alpha \cap U_\beta$ and $v'
= \rho_{\alpha, \beta}(p)(v)$, where $\rho_{\alpha, \beta}(p)(v)$ is
the value of the linear transformation on ${\bf R}$ associated to
$\rho_{\alpha, \beta}(p)$ at $v$.  This leads to a real vector bundle
$E$ of rank $n$ over $X$.

	Suppose that for each $\alpha \in A$ there is a continuous
mapping $\tau_\alpha$ from $U_\alpha$ into $GL(n, {\bf R})$ such that
$\rho_{\alpha, \beta} = \tau_\beta \, \tau_\alpha^{-1}$ on $U_\alpha
\cap U_\beta$ when this set is nonempty.  If $p \in U_\alpha \cap
U_\beta$ and $v \in {\bf R}^n$, then $(p, \tau_\alpha(p)(v)) \in
E_\alpha$ and $(p, \tau_\beta(p)(v)) \in E_\beta$ are identified in
$E$.  Hence we get an isomorphism from $X \times {\bf R}^n$ onto $E$
by sending $(p, v) \in X \times {\bf R}^n$ to $(p, \tau_\alpha(p)(v))
\in E_\alpha$ when $p \in U_\alpha$.

	Similar remarks apply to complex matrices and vector bundles.
One can also view this as a special case of the previous situation, by
identifying $n \times n$ complex matrices with $2n \times 2n$ real
matrices.  Other groups of matrices correspond to other types of
vector bundles.  A vector bundle is said to be \emph{flat} if the
transition functions $\rho_{\alpha, \beta}$ can be chosen to be
constant.

\section{Subgroups}
\label{subgroups}
\setcounter{equation}{0}

	Let $G$ be a group, and let $H$ be a subgroup of $G$.  The
quotient space $G / H$ consists of the cosets $a H = \{ a \, b : b \in
H \}$ of $H$ in $G$ for $a \in G$, and the quotient mapping $q : G \to
H$ sends $a \in G$ to $a H$.  If $G$ is a topological group and $H$ is
a closed subgroup of $G$, then $G / H$ can be given the associated
quotient topology, for which the quotient mapping is continuous and
open.

	Let $X$ be a topological space, and let $\mathcal{U} =
\{U_\alpha\}_{\alpha \in A}$ be an open covering of $X$.  Suppose that
$\phi$ is a continuous mapping from $X$ into $G / H$.  Under suitable
conditions, for each $\alpha \in A$ there is a continuous mapping
$\phi_\alpha$ from $U_\alpha$ into $G$ whose composition with the
quotient mapping $q$ is equal to $\phi$ on $U_\alpha$.  If $\alpha,
\beta \in A$ and $U_\alpha \cap U_\beta \ne \emptyset$, then
\begin{equation}
\label{rho_{alpha, beta} = phi_alpha^{-1} phi_beta}
	\rho_{\alpha, \beta} = \phi_\alpha^{-1} \, \phi_\beta
\end{equation}
is a continuous mapping from $U_\alpha \cap U_\beta$ into $G$ such
that $\rho_{\beta, \alpha} = \rho_{\alpha, \beta}^{-1}$ and
$\rho_{\alpha, \beta} \, \rho_{\beta, \gamma} \, \rho_{\gamma, \alpha}
= e$ on $U_\alpha \cap U_\beta \cap U_\gamma$ when this set is
nonempty.  Moreover, $\rho_{\alpha, \beta}(p) \in H$ for every $p \in
U_\alpha \cap U_\beta$, since $q(\phi_\alpha(p)) = q(\phi_\beta(p)) =
\phi(p)$.

	Conversely, suppose that $\phi_\alpha : U_\alpha \to G$ is a
continuous mapping for each $\alpha \in A$, and let $\rho_{\alpha,
\beta}$ be as in (\ref{rho_{alpha, beta} = phi_alpha^{-1} phi_beta}).
As before, $\rho_{\beta, \alpha} = \rho_{\alpha, \beta}^{-1}$ and
$\rho_{\alpha, \beta} \, \rho_{\beta, \gamma} \, \rho_{\gamma, \alpha}
= e$ on $U_\alpha \cap U_\beta \cap U_\gamma$ when this set is
nonempty.  If every $\rho_{\alpha, \beta}$ takes values in $H$, then
$q \circ \phi_\alpha = q \circ \phi_\beta$ on $U_\alpha \cap U_\beta$
and $q(\phi_\alpha(p))$ determines a continuous mapping from $X$ into
$G / H$.  If for each $\alpha \in A$ there is also a continuous
mapping $\tau_\alpha : U_\alpha \to H$ such that $\rho_{\alpha, \beta}
= \tau_\alpha^{-1} \, \tau_\beta$ on $U_\alpha \cap U_\beta$ when this
set is nonempty, and if we put $\phi'_\alpha = \phi_\alpha \,
\tau_\alpha^{-1}$, then $q \circ \phi'_\alpha = q \circ \phi_\alpha$,
$\phi'_\alpha = \phi'_\beta$ on $U_\alpha \cap U_\beta$ when this set
is nonempty, and $\phi'_\alpha(p)$ determines a continuous mapping
from $X$ into $G$.


\begin{thebibliography}{74}


\bibitem {a-k} L.~Ambrosio and B.~Kirchheim, {\it Currents in metric
spaces}, Acta Mathematica {\bf 185} (2000), 1--80.

\bibitem {a-t} L.~Ambrosio and P.~Tilli, {\it Topics in Analysis in
Metric Spaces}, Oxford University Press, 2004.

\bibitem {a} M.~Atiyah, {\it $K$-Theory}, based on notes by
D.~Anderson, 2nd edition, Addison-Wesley, 1989.

\bibitem {au} T.~Aubin, {\it A Course in Differential Geometry},
American Mathematical Society, 2001.

\bibitem {a-c-g} P.~Auscher, T.~Coulhon, and A.~Grigoryan, editors,
{\it Heat Kernels and Analysis on Manifolds, Graphs, and Metric
Spaces}, American Mathematical Society, 2003.

\bibitem {b} M.~Berger, {\it Geometry I, II}, translated from the
French by M.~Cole and S.~Levy, Springer-Verlag, 1987.

\bibitem {b'} M.~Berger, {\it A Panoramic View of Riemannian
Geometry}, Springer-Verlag, 2003.

\bibitem {b-g} M.~Berger and B.~Gostiaux, {\it Differential Geometry:
Manifolds, Curves, and Surfaces}, translated from the French by
S.~Levy, Springer-Verlag, 1988.

\bibitem {bth} W.~Boothby, {\it An Introduction to Differentiable
Manifolds and Riemannian Geometry}, 2nd edition, Academic Press, 1986.

\bibitem {b-t} R.~Bott and L.~Tu, {\it Differential Forms in Algebraic
Topology}, Springer-Verlag, 1982.

\bibitem {br1} G.~Bredon, {\it Sheaf Theory}, 2nd edition,
Springer-Verlag, 1997.

\bibitem {br2} G.~Bredon, {\it Topology and Geometry},
Springer-Verlag, 1997.

\bibitem {c1} M.~do Carmo, {\it Differential Geometry of Curves and
Surfaces}, translated from the Portuguese, Prentice-Hall, 1976.

\bibitem {c2} M.~do Carmo, {\it Riemannian Geometry}, translated from
the second Portuguese edition by F.~Flaherty, Birkh\"auser, 1992.

\bibitem {c3} M.~do Carmo, {\it Differential Forms and Applications},
translated from the 1971 Portuguese original, Springer-Verlag, 1994.

\bibitem {c-w-1} R.~Coifman and G.~Weiss, {\it Analyse Harmonique
Non-Commutative sur Certains Espaces Homog\`enes}, Lecture Notes in
Mathematics {\bf 242}, Springer-Verlag, 1971.

\bibitem {c-w-2} R.~Coifman and G.~Weiss, {\it Extensions of Hardy
spaces and their use in analysis}, Bulletin of the American
Mathematical Society {\bf 83} (1977), 569--645.

\bibitem {cnl} L.~Conlon, {\it Differentiable Manifolds}, 2nd edition,
Birkh\"auser, 2001.

\bibitem {cns} A.~Connes, {\it Noncommutative Geometry}, Academic
Press, 1994.

\bibitem {fa} K.~Falconer, {\it The Geometry of Fractal Sets},
Cambridge University Press, 1986.

\bibitem {fe} H.~Federer, {\it Geometric Measure Theory},
Springer-Verlag, 1969.

\bibitem {fl} H.~Flanders, {\it Differential Forms with Applications
to the Physical Sciences}, Dover, 1989.

\bibitem {g-g} T.~Gamelin and R.~Greene, {\it Introduction to
Topology}, 2nd edition, Dover, 1999.

\bibitem {g-j} L.~Gillman and M.~Jerison, {\it Rings of Continuous
Functions}, Springer-Verlag, 1976.

\bibitem {gi} E.~Giusti, {\it Minimal Surfaces and Functions of
Bounded Variation}, Birkh\"auser, 1984.

\bibitem {gra-r} H.~Grauert and R.~Remmert, {\it Coherent Analytic
Sheaves}, Springer-Verlag, 1984.

\bibitem {gry} B.~Gray, {\it Homotopy Theory: An Introduction to
Algebraic Topology}, Academic Press, 1975.

\bibitem {g} R.~Gunning, {\it Introduction to Holomorphic Functions of
Several Variables}, Volumes 1, 2, 3, Wadsworth \& Brooks / Cole, 1990.

\bibitem {g-r} R.~Gunning and H.~Rossi, {\it Analytic Functions of
Several Complex Variables}, Prentice-Hall, 1965.

\bibitem {ht} A.~Hatcher, {\it Algebraic Topology}, Cambridge
University Press, 2002.

\bibitem {h1} J.~Heinonen, {\it Lectures on Analysis on Metric Spaces},
Springer-Verlag, 2001.

\bibitem {h2} J.~Heinonen, {\it Geometric Embeddings of Metric
Spaces}, Reports of the Department of Mathematics and Statistics {\bf
90}, University of Jyv\"askyl\"a, 2003.

\bibitem {h3} J.~Heinonen, {\it Lectures on Lipschitz Analysis},
Reports of the Department of Mathematics and Statistics {\bf 100},
University of Jyv\"askyl\"a, 2005.

\bibitem {hrs} M.~Hirsch, {\it Differential Topology},
Springer-Verlag, 1994.

\bibitem {hrm} L.~H\"ormander, {\it An Introduction to Complex
Analysis in Several Complex Variables}, 3rd edition, North-Holland,
1990.

\bibitem {h-w} W.~Hurewicz and H.~Wallman, {\it Dimension Theory},
revised edition, Princeton University Press, 1948.

\bibitem {kap} I.~Kaplansky, {\it Set Theory and Metric Spaces}, 2nd
edition, Chelsea, 1977.

\bibitem {k-n} S.~Kobayashi and K.~Nomizu, {\it Foundations of
Differential Geometry}, Volumes I and II, Wiley, 1996.

\bibitem {k} S.~Krantz, {\it Function Theory of Several Complex
Variables}, AMS Chelsea, 2001.

\bibitem {k-p-1} S.~Krantz and H.~Parks, {\it The Geometry of Domains in
Space}, Birkh\"auser, 1999.

\bibitem {k-p-2} S.~Krantz and H.~Parks, {\it The Implicit Function
Theorem: History, Theory, and Applications}, Birkh\"auser, 2002.

\bibitem {l} S.~Lang, {\it Introduction to Differentiable Manifolds},
2nd edition, Springer-Verlag, 2002.

\bibitem {l1} J.~Lee, {\it Riemannian Manifolds: An Introduction to
Curvature}, Springer-Verlag, 1997.

\bibitem {l2} J.~Lee, {\it Introduction to Topological Manifolds},
Springer-Verlag, 2000.

\bibitem {l3} J.~Lee, {\it Introduction to Smooth Manifolds},
Springer-Verlag, 2003.

\bibitem {l-v} J.~Luukkainen and J.~V\"ais\"al\"a, {\it Elements of
Lipschitz topology}, Annales Academiae Scientiarum Fennicae Series A I
Mathematica {\bf 3} (1977), 85--122.

\bibitem {mac} S.~MacLane, {\it Categories for the Working
Mathematician}, 2nd edition, Springer-Verlag, 1998.

\bibitem {m-t} I.~Madsen and J.~Tornehave, {\it From Calculus to
Cohomology: de Rham Cohomology and Characteristic Classes}, Cambridge
University Press, 1997.

\bibitem {ms1} W.~Massey, {\it Homology and Cohomology Theory: An
Approach Based on Alexander--Spanier Cochains}, Dekker, 1978.

\bibitem {ms2} W.~Massey, {\it A Basic Course in Algebraic Topology},
Springer-Verlag, 1991.

\bibitem {mt} P.~Mattila, {\it Geometry of Sets and Measures in
Euclidean Spaces: Fractals and Rectifiability}, Cambridge University
Press, 1995.

\bibitem {mc} G.~McCarty, {\it Topology: An Introduction with
Application to Topological Groups}, 2nd edition, Dover, 1988.

\bibitem {mn} B.~Mendelson, {\it Introduction to Topology}, 3rd
edition, Dover, 1990.

\bibitem {mln} J.~Milnor, {\it Topology from the Differentiable
Viewpoint}, based on notes by D.~Weaver, Princeton University Press,
1997.

\bibitem {m-s} J.~Milnor and J.~Stasheff, {\it Characteristic
Classes}, Princeton University Press, 1974.

\bibitem {mo1} F.~Morgan, {\it Riemannian Geometry: A Beginner's
Guide}, 2nd edition, A K Peters, 1998.

\bibitem {mo2} F.~Morgan, {\it Geometric Measure Theory: A Beginner's
Guide}, 3rd edition, Academic Press, 2000.

\bibitem {mk} A.~Mukherjee, {\it Topics in Differential Topology},
Hindustan Book Agency, 2005.

\bibitem {p} A.~Pressley, {\it Elementary Differential Geometry},
Springer-Verlag, 2001.

\bibitem {rh} G.~de Rham, {\it Differentiable Manifolds: Forms,
Currents, Harmonic Forms}, translated from the French by F.~Smith,
with an introduction by S.~Chern, Springer-Verlag, 1984.

\bibitem {r} W.~Rudin, {\it Principles of Mathematical Analysis}, 3rd
edition, McGraw-Hill, 1976.

\bibitem {si} G.~Simmons, {\it Introduction to Topology and Modern
Analysis}, Krieger, 1983.

\bibitem {sp} E.~Spanier, {\it Algebraic Topology}, Springer-Verlag,
1981.

\bibitem {spv} M.~Spivak, {\it Calculus on Manifolds: A Modern
Approach to Classical Theorems of Advanced Calculus}, Benjamin, 1965.

\bibitem {s} N.~Steenrod, {\it The Topology of Fiber Bundles},
Princeton University Press, 1999.

\bibitem {s-e} N.~Steenrod {\it Cohomology Operations}, based on
lectures by N.~Steenrod and written and revised by D.~Epstein,
Princeton University Press, 1962.

\bibitem {st1} E.~Stein, {\it Singular Integrals and Differentiability
Properties of Functions}, Princeton University Press, 1970.

\bibitem {st2} E.~Stein, {\it Harmonic Analysis: Real-Variable
Methods, Orthogonality, and Oscillatory Integrals}, with the
assistance of T.~Murphy, Princeton University Press, 1993.

\bibitem {s-w} E.~Stein and G.~Weiss, {\it Introduction to Fourier
Analysis on Euclidean Spaces}, Princeton University Press, 1971.

\bibitem {v} J.~V\"ais\"al\"a, {\it Metric duality in Euclidean
spaces}, Mathematica Scandinavica {\bf 80} (1997), 249--288.

\bibitem {w} F.~Warner, {\it Foundations of Differential Manifolds and
Lie Groups}, Springer-Verlag, 1983.

\bibitem {wells} R.~Wells, {\it Differential Analysis on Complex
Manifolds}, 2nd edition, Springer-Verlag, 1980.

\bibitem {wh} G.~Whitehead, {\it Elements of Homotopy Theory},
Springer-Verlag, 1978.

\bibitem {wn} H.~Whitney, {\it Geometric Integration Theory},
Princeton University Press, 1957.




\end{thebibliography}
\end{document}